\documentclass[12pt,a4paper]{article}

\usepackage{a4wide}
\usepackage{amsfonts,amsmath}
\usepackage{latexsym}

\usepackage{xcolor}

\usepackage{times}
\author{Karoly J. Boroczk\footnote{Alfred Renyi Institute of Mathematics, Realtanoda utca 13-15, 1053, Budapest, Hungary, Supported by NKFIH 132002}, Pavlos Kalantzopoulos\footnote{Central European University, Nador utca 9, 1051, Budapest, Hungary}, Dongmeng Xi\footnote{Department of Mathematics,
Shanghai University,
Shanghai 200444, China, Supported by National Natural Science Foundation of China (12071277).}}
\title{The case of equality in geometric instances of Barthe's reverse Brascamp-Lieb inequality}

\newcommand{\proof}{\noindent{\it Proof: }}
\newcommand{\proofbox}{\mbox{ $\Box$}\\}
\newcommand{\R}{\mathbb{R}}

\newtheorem{lemma}{Lemma}
\newtheorem{theo}[lemma]{Theorem}
\newtheorem{example}[lemma]{Example}

\newtheorem{coro}[lemma]{Corollary}

\newtheorem{prop}[lemma]{Proposition}

\begin{document}

\maketitle

\begin{abstract}
 The works of Bennett, Carbery, Christ, Tao and of Valdimarsson have clarified when equality holds in the Brascamp-Lieb inequality.
Here we characterize the case of equality in the Geometric case of Barthe's reverse Brascamp-Lieb inequality.
\end{abstract}

\section{Introduction}

For a proper linear subspace $E$ of $\R^n$ ($E\neq \R^n$ and $E\neq\{0\}$),
let $P_E$ denote the orthogonal projection into $E$.
We say that the  subspaces $E_1,\ldots,E_k$ of $\R^n$ and $c_1,\ldots,c_k>0$  form
a Geometric Brascamp-Lieb data if they satisfy
\begin{equation}
\label{highdimcond0}
\sum_{i=1}^kc_iP_{E_i}=I_n.
\end{equation}
The name ``Geometric Brascamp-Lieb data" coined by Bennett, Carbery, Christ, Tao \cite{BCCT08} comes from the following theorem, originating in the work of Brascamp, Lieb \cite{BrL76} and Ball \cite{Bal89,Bal91} in the rank one case
(${\rm dim}\,E_i=1$ for $i=1,\ldots,k$),
and Lieb \cite{Lie90} and Barthe \cite{Bar98} in the general case. In the rank one case, the Geometric Brascamp-Lieb data is known as Parseval frame in coding theory and computer science (see for example Casazza, Tran, Tremain \cite{CTT20}).

\begin{theo}[Brascamp-Lieb, Ball, Barthe]
\label{BLtheo}
For  the linear subspaces  $E_1,\ldots,E_k$ of $\R^n$ and $c_1,\ldots,c_k>0$ satisfying
\eqref{highdimcond0}, and for non-negative $f_i\in L_1(E_i)$, we have
\begin{equation}
\label{BL}
\int_{\R^n}\prod_{i=1}^kf_i(P_{E_i}x)^{c_i}\,dx
\leq \prod_{i=1}^k\left(\int_{E_i}f_i\right)^{c_i}
\end{equation}
\end{theo}
{\bf Remark} This is H\"older's inequality if $E_1=\ldots=E_k=\R^n$ and $B_i=I_n$, and hence
$\sum_{i=1}^kc_i=1$.\\

We note that equality holds in Theorem~\ref{BLtheo}  if $f_i(x)=e^{-\pi\|x\|^2}$ for $i=1,\ldots,k$; and hence, each $f_i$ is a Gaussian density.
Actually,
Theorem~\ref{BLtheo} is an important special case discovered by Ball \cite{Bal91,Bal03}
in the rank one case and by Barthe \cite{Bar98} in the general case of the general Brascamp-Lieb inequality Theorem~\ref{BLgeneral}.

After partial results by Barthe \cite{Bar98}, Carlen,  Lieb,  Loss \cite{CLL04}
and Bennett, Carbery, Christ, Tao \cite{BCCT08}, it was Valdimarsson \cite{Val08}
who characterized equality in the Geometric Brascamp-Lieb inequality. In order to state his result, we need some notation.
Let $E_1,\ldots,E_k$  the proper linear subspaces of $\R^n$  and $c_1,\ldots,c_k>0$ satisfy
\eqref{highdimcond0}.
In order to understand extremizers in \eqref{RBL}, following
Carlen,  Lieb,  Loss \cite{CLL04} and Bennett, Carbery, Christ, Tao \cite{BCCT08}, we say that a non-zero linear subspace $V$
is a
critical subspace  if
$$
\sum_{i=1}^k c_i\dim(E_i\cap V)=\dim V,
$$
which is turn equivalent saying that
$$
\mbox{$E_i=(E_i\cap V)+ (E_i\cap V^\bot)$
for $i=1,\ldots,k$}
$$
according to \cite{BCCT08} (see also Lemma~\ref{highdimcritical}). We say that a critical subspace $V$ is indecomposable if
$V$ has no proper critical linear subspace.

Valdimarsson \cite{Val08} introduced the so called independent subspaces and the dependent space.
We write $J$ to denote the set of $2^k$ functions
$\{1,\ldots,k\}\to\{0,1\}$. If $\varepsilon\in J$, then
let $F_{(\varepsilon)}=\cap_{i=1}^kE_i^{(\varepsilon(i))}$ where $E_i^{(0)}=E_i$ and $E_i^{(1)}=E_i^\bot$
for $i=1,\ldots,k$. We write $J_0$ to denote the subset of $\varepsilon\in J$ such that
${\rm dim}\,F_{(\varepsilon)}\geq 1$, and such an $F_{(\varepsilon)}$  is called independent following Valdimarsson \cite{Val08}. Readily $F_{(\varepsilon)}$ and $F_{(\tilde{\varepsilon})}$ are orthogonal
if $\varepsilon\neq\tilde{\varepsilon}$ for $\varepsilon,\tilde{\varepsilon}\in J_0$.
 In addition, we write $F_{\rm dep}$ to denote the orthogonal component of
$\oplus_{\varepsilon \in J_0}F_{(\varepsilon)}$. In particular, $\R^n$ can be written as a direct sum of pairwise orthogonal linear subspaces in the form
\begin{equation}
\label{independent-dependent0}
\R^n=\left(\oplus_{\varepsilon \in J_0}F_{(\varepsilon)}\right)\oplus F_{\rm dep}.
\end{equation}
Here it is possible that $J_0=\emptyset$, and hence $\R^n=F_{\rm dep}$, or
$F_{\rm dep}=\{0\}$, and hence $\R^n=\oplus_{\varepsilon \in J_0}F_{(\varepsilon)}$ in that case.

For a non-zero linear subspace $L\subset \R^n$, we say that a linear transformation $A:\,L\to L$ is positive definite
if  $\langle Ax,y\rangle=\langle x, Ay\rangle$ and $\langle x, Ax\rangle>0$ for any $x,y\in L\backslash\{0\}$.

\begin{theo}[Valdimarsson]
\label{BLtheoequa}
For  the proper linear subspaces  $E_1,\ldots,E_k$ of $\R^n$ and $c_1,\ldots,c_k>0$ satisfying
\eqref{highdimcond0}, let us assume that equality holds in the Brascamp-Lieb inequality \eqref{BL}
for non-negative $f_i\in L_1(E_i)$, $i=1,\ldots,k$.
If $F_{\rm dep}\neq\R^n$, then let $F_1,\ldots,F_\ell$ be the independent subspaces, and if
$F_{\rm dep}=\R^n$, then let $\ell=1$ and $F_1=\{0\}$. There exist
$b\in F_{\rm dep}$ and $\theta_i>0$ for $i=1,\ldots,k$,
integrable non-negative $h_{j}:\,F_j\to[0,\infty)$  for $j=1,\ldots,\ell$, and a positive definite matrix
$A:F_{\rm dep}\to F_{\rm dep}$ such that
the eigenspaces of $A$ are critical subspaces and
\begin{equation}
\label{BLtheoequaform}
f_i(x)=\theta_i e^{-\langle AP_{F_{\rm dep}}x,P_{F_{\rm dep}}x-b\rangle}\prod_{F_j\subset E_i}h_{j}(P_{F_j}(x))
\mbox{ \ \ \  for Lebesgue a.e. $x\in E_i$}.
\end{equation}
On the other hand, if for any $i=1,\ldots,k$, $f_i$ is of the form as in \eqref{BLtheoequaform}, then equality holds in \eqref{BL} for $f_1,\ldots,f_k$.
\end{theo}

Theorem~\ref{BLtheoequa} explains the term "independent subspaces" because the functions
$h_{j}$ on $F_j$ are chosen freely and independently  from each other.

A reverse form of the Geometric Brascamp-Lieb inequality was proved by
Barthe \cite{Bar98}. We write $\int_{*,\R^n}\varphi $ to denote the
inner integral for a possibly non-integrable function $\varphi:\,\R^n\to[0,\infty)$; namely, the supremum (actually maximum)
of $\int_{\R^n} \psi$ where $0\leq \psi\leq \varphi$ is Lebesgue measurable.

\begin{theo}[Barthe]
\label{RBLtheo}
For  the non-trivial linear subspaces  $E_1,\ldots,E_k$ of $\R^n$ and $c_1,\ldots,c_k>0$ satisfying
\eqref{highdimcond0}, and for non-negative $f_i\in L_1(E_i)$, we have
\begin{equation}
\label{RBL}
\int_{*,\R^n}\sup_{x=\sum_{i=1}^kc_ix_i,\, x_i\in E_i}\;\prod_{i=1}^kf_i(x_i)^{c_i}\,dx
\geq \prod_{i=1}^k\left(\int_{E_i}f_i\right)^{c_i}.
\end{equation}
\end{theo}
{\bf Remark} This is the Pr\'ekopa-Leindler inequality Theorem~\ref{PL}
 if $E_1=\ldots=E_k=\R^n$ and $B_i=I_n$, and hence $\sum_{i=1}^kc_i=1$.\\

The function $x\mapsto \sup_{x=\sum_{i=1}^kc_ix_i,\, x_i\in E_i}\;\prod_{i=1}^kf_i(x_i)^{c_i}$ in Theorem~\ref{RBLtheo} 
 may not be measurable if
$f_1,\ldots,f_k\geq 0$ are measurable but it is measurable if $f_1,\ldots,f_k$ are Borel (and hence analytic in the set theoretic sense).
We note that Theorem~\ref{RBLtheo} is usually stated using outer integrals (even in  \cite{Bar98}) but it actually holds for inner integrals if one closely follows Barthe's argument. 
If $f_1,\ldots,f_k$ are $C^1$ and positive, then \cite{Bar98} provides an elegant argument using optimal transport. 
For general measurable $f_1,\ldots,f_k$, one can approximate each $f_i$ by non-negative linear combination of characteristic functions of compact sets from below, and hence enoungh to consider such functions in order to prove Theorem~\ref{RBLtheo}.
However, the characteristic function of a compact set can be approximated from above by a decreasing sequence of positive $C^\infty$ functions, completing the proof of  Theorem~\ref{RBLtheo}.

 We say that a function $h:\,\R^n\to[0,\infty)$
is log-concave if $h((1-\lambda)x+\lambda\,y)\geq h(x)^{1-\lambda}h(y)^\lambda$ for any
$x,y\in\R^n$ and $\lambda\in(0,1)$; or in other words, $h=e^{-W}$ for a convex function $W:\,\R^n\to(-\infty,\infty]$.
Our main result is the following characterization of equality in the Geometric Barthe's inequality
\eqref{RBL}.

\begin{theo}
\label{RBLtheoequa}
For linear subspaces  $E_1,\ldots,E_k$ of $\R^n$ and $c_1,\ldots,c_k>0$ satisfying
\eqref{highdimcond0},
if $F_{\rm dep}\neq\R^n$, then let $F_1,\ldots,F_\ell$ be the independent subspaces, and if
$F_{\rm dep}=\R^n$, then let $\ell=1$ and $F_1=\{0\}$.

 If equality holds in the Geometric Barthe's inequality \eqref{RBL}
for non-negative $f_i\in L_1(E_i)$ with $\int_{E_i}f_i>0$, $i=1,\ldots,k$, then
\begin{equation}
\label{RBLtheoequaform}
f_i(x)=\theta_i e^{-\langle AP_{F_{\rm dep}}x,P_{F_{\rm dep}}x-b_i\rangle}\prod_{F_j\subset E_i}h_{j}(P_{F_j}(x-w_i)) \mbox{ \ \ \  for Lebesgue a.e. $x\in E_i$}
\end{equation}
where
\begin{itemize}
\item $\theta_i>0$,
$b_i\in E_i\cap F_{\rm dep}$ and $w_i\in E_i$ for $i=1,\ldots,k$,
\item $h_{j}\in L_1(F_j)$ is non-negative for $j=1,\ldots,\ell$, and in addition,
$h_j$ is log-concave if there exist $\alpha\neq \beta$ with $F_j\subset E_\alpha\cap E_\beta$,
\item $A:F_{\rm dep}\to F_{\rm dep}$ is a positive definite matrix such that
the eigenspaces of $A$ are critical subspaces.
\end{itemize}
On the other hand, if for any $i=1,\ldots,k$, $f_i$ is of the form as in \eqref{RBLtheoequaform} and equality holds for  all
$x\in E_i$ in \eqref{RBLtheoequaform}, then equality holds in \eqref{RBL} for $f_1,\ldots,f_k$.
\end{theo}
{\bf Remark} An independent subspace $F_j$ is contained in a single $E_\alpha$ if and only if $E_i\subset E_\alpha^\bot$ for $i\neq\alpha$, and hence $F_j=E_\alpha$. In particular, if for any $\alpha=1,\ldots,k$, $\{E_i\}_{i\neq \alpha}$ spans $\R^n$ in Theorem~\ref{RBLtheoequa},
then any extremizer of the Geometric Barthe's inequality is log-concave.\\

The explanation for the phenomenon concerning the log-concavity of $h_j$ in Theorem~\ref{RBLtheoequa} is as follows (see the proof of Proposition~\ref{Thetasplitshij}). Let
$\ell\geq 1$ and $j\in\{1,\ldots,\ell\}$, and hence $\sum_{E_i\supset F_j}c_i=1$. If $f_1,\ldots,f_k$ are of the form \eqref{RBLtheoequaform}, then equality in  Barthe's inequality \eqref{RBL} yields
$$
\int_{*,F_j}\sup_{x=\sum_{E_i\supset F_j}c_i x_i\atop x_i\in F_j}h_{j}\Big(x_i-P_{F_j}w_i\Big)^{c_i}\,dx=
\prod_{E_i\supset F_j}\left(\int_{F_j}h_{j}\Big(x-P_{F_j}w_i\Big)\,dx\right)^{c_i}
\left(= \int_{F_j} h_j(x)\,dx\right).
$$
Therefore, if there exist $\alpha\neq \beta$ with $F_j\subset E_\alpha\cap E_\beta$, then
the equality conditions in the Pr\'ekopa-Leindler inequality
Proposition~\ref{PL} imply that $h_j$ is log-concave. On the other hand, if there exists $\alpha\in \{1,\ldots,k\}$ such that
$F_j\subset E_\beta^\bot$ for   $\beta\neq\alpha$, then we do not have any condition on $h_j$, and $c_\alpha=1$.\\

For completeness, let us state and discuss the general Brascamp-Lieb inequality and its reverse form due to Barthe.
The following was proved by Brascamp, Lieb \cite{BrL76}  in the rank one case and Lieb \cite{Lie90} in general.

\begin{theo}[Brascamp-Lieb Inequality]
\label{BLgeneral}
Let $B_i:\R^n\to H_i$ be surjective linear maps where $H_i$ is $n_i$-dimensional Euclidean space,
$n_i\geq 1$, for $i=1,\ldots,k$,
and let $c_1,\ldots,c_k>0$ satisfy $\sum_{i=1}^kc_in_i=n$.
For non-negative $f_i\in L_1(H_i)$, we have
\begin{equation}
\label{BLgeneraleq}
\int_{\R^n}\prod_{i=1}^kf_i(B_ix)^{c_i}\,dx
\leq C\prod_{i=1}^k\left(\int_{H_i}f_i\right)^{c_i}
\end{equation}
where $C$ is determined by choosing centered Gaussians $f_i(x)=e^{-\langle A_ix,x\rangle}$, $A_i$ positive definite.
\end{theo}
{\bf Remark} The Geometric Brascamp-Lieb Inequality is readily a special case of
\eqref{BLgeneraleq}. We note that \eqref{BLgeneraleq}
 isH\"older's inequality if $H_1=\ldots=H_k=\R^n$ and each $B_i=I_n$, and hence $C=1$
and $\sum_{i=1}^kc_i=1$ in that case.\\

We say that two Brascamp-Lieb data $\{(B_i,c_i)\}_{i=1,\ldots,k}$ and
$\{(B'_i,c'_i)\}_{i=1,\ldots,k'}$ as in Theorem~\ref{BLgeneral}
are called equivalent if $k'=k$, $c'_i=c_i$, and there exists linear isomorphism $\Phi_i:H_i\to H'_i$ for
$i=1,\ldots,k$ such that $B'_i=\Phi_i\circ B_i$. It was proved by
Carlen,  Lieb,  Loss \cite{CLL04} in the rank one case, and by Bennett, Carbery, Christ, Tao \cite{BCCT08}
in general that there exists a set of extremizers $f_1,\ldots,f_k$ for \eqref{BLgeneraleq} if and only
if the Brascamp-Lieb data $\{(B_i,c_i)\}_{i=1,\ldots,k}$ is equivalent to some Geometric Brascamp-Lieb data.
Therefore, Valdimarsson's Theorem~\ref{BLtheoequa} provides a full characterization of the equality case
in Theorem~\ref{BLgeneral}, as well.

The following reverse version of the Brascamp-Lieb inequality was proved by Barthe in \cite{Bar97} in the rank one case, and
in \cite{Bar98} in general.

\begin{theo}[Barthe's Inequality]
\label{RBLgeneral}
Let $B_i:\R^n\to H_i$ be surjective linear maps where $H_i$ is $n_i$-dimensional Euclidean space,
$n_i\geq 1$, for $i=1,\ldots,k$,
and let $c_1,\ldots,c_k>0$ satisfy $\sum_{i=1}^kc_in_i=n$.
For non-negative $f_i\in L_1(H_i)$, we have
\begin{equation}
\label{RBLgeneraleq}
\int_{*,\R^n}
\sup_{x=\sum_{i=1}^kc_i B_i^*x_i,\, x_i\in H_i}\;
\prod_{i=1}^kf_i(x_i)^{c_i}\,dx
\geq D\prod_{i=1}^k\left(\int_{H_i}f_i\right)^{c_i}
\end{equation}
where $D$ is determined by choosing centered Gaussians $f_i(x)=e^{-\langle A_ix,x\rangle}$, $A_i$ positive definite.
\end{theo}
{\bf Remark} The Geometric  Barthe's Inequality \eqref{RBL} is readily a special case of
\eqref{RBLgeneraleq}. We note that \eqref{RBLgeneraleq} is the Pr\'ekopa-Leindler inequality if $H_1=\ldots=H_k=\R^n$
 and each $B_i=I_n$, and hence $D=1$ and $\sum_{i=1}^kc_i=1$ in that case.\\

Actually, again, Barthe  \cite{Bar98} stated \eqref{RBLgeneraleq} with outer integrals, but his argument yields the same stament with inner inner integrals, see the discussion after Theorem~\ref{RBLtheo}.

Concerning extremals in Theorem~\ref{RBLgeneral},
Lehec \cite{Leh14} proved that if there exists some Gaussian extremizers for
Barthe's Inequality \eqref{RBLgeneraleq}, then the corresponding Brascamp-Lieb data
 $\{(B_i,c_i)\}_{i=1,\ldots,k}$ is equivalent to some Geometric Brascamp-Lieb data; therefore, the equality case of
\eqref{RBLgeneraleq} can be understood via Theorem~\ref{RBLtheoequa} in that case.

However, it is still not known whether having any extremizers in Barthe's Inequality \eqref{RBLgeneraleq} yields the existence of Gaussian extremizers. One possible approach is to use iterated convolutions and
renormalizations as in Bennett, Carbery, Christ, Tao \cite{BCCT08} in the case of Brascamp-Lieb inequality.

There are three main methods of proofs that work for proving both
the Brascamp-Lieb Inequality and its reverse form, Barthe's inequality. The paper Barthe \cite{Bar98} used optimal transportation to prove
Barthe's Inequality (``the Reverse Brascamp-Lieb inequality") and reprove the Brascamp-Lieb Inequality simultaneously.
A heat equation argument was provided in the rank one case by
Carlen,  Lieb,  Loss \cite{CLL04} for the Brascamp-Lieb Inequality and by
Barthe, Cordero-Erausquin \cite{BaC04} for Barthe's inequality. The general versions of both inequalities are proved via the heat equation approach by
Barthe, Huet \cite{BaH09}. Finally, simultaneous probabilistic arguments for the two inequalities are due to
Lehec \cite{Leh14}.

We note that
Chen,  Dafnis, Paouris \cite{CDP15}
and Courtade, Liu \cite{CoL21}, as well, deal systematically with finiteness conditions in Brascamp-Lieb and
Barthe's inequalities. The importance of the Brascamp-Lieb inequality is shown by the fact that besides harmonic analysis, probability and convex geometry,
it has been also even applied in number theory, see eg. Guo, Zhang \cite{GuZ19}.
Various versions of the Brascamp-Lieb inequality and its reverse form have been obtained by
Balogh, Kristaly \cite{BaK18}
Barthe \cite{Bar04}, Barthe, Cordero-Erausquin \cite{BaC04},
Barthe, Cordero-Erausquin,  Ledoux, Maurey \cite{BCLM11},
Barthe, Wolff \cite{BaW14,BaW22},
Bennett, Bez, Flock, Lee \cite{BBFL18},
 Bennett, Bez, Buschenhenke, Cowling, Flock \cite{BBBCF20},
Bobkov, Colesanti, Fragal\`a \cite{BCF14},
Bueno, Pivarov \cite{BuP21},
Chen,  Dafnis, Paouris \cite{CDP15},  Courtade, Liu \cite{CoL21},
Duncan \cite{Dun21}, Ghilli, Salani \cite{GhS17},
Kolesnikov, Milman \cite{KoM},  Livshyts \cite{Liv21,Liv},
Lutwak, Yang, Zhang \cite{LYZ04,LYZ07},
Maldague \cite{Mal},  Marsiglietti \cite{Mar17}, Rossi, Salani \cite{RoS17,RoS19}.

Concerning the proof of Theorem~\ref{RBLtheoequa},
we discuss the structure theory of a Brascamp-Lieb data, Barthe's crucial determinantal inequality
({\it cf.} Proposition~\ref{GBLconstanthighdim})
and the extremality of Gaussians ({\it cf.} Proposition~\ref{GaussianExtremizers})
in Sections~\ref{secStructure}, \ref{secStuctureHigh} and \ref{secGaussians}.
Section~\ref{secSplitting} explains how Barthe's proof of his inequality using optimal transportation in \cite{Bar98}
yields the splitting along independent and dependent subspaces in the case of equality in Barthe's inequality
for positive $C^1$ probality densities $f_1,\ldots,f_k$, and how the equality case of the Pr\'ekopa-Leindler inequality leads to the log-concavity of certain functions involved. However, one still needs to produce suitably smooth extremizers given any extremizers of Barthe's inequality. In order to achieve this, we discuss that convolution and suitable products of extremizers are also extremizers in Section~\ref{secExtremizers}. To show that extremizers are Gaussians on the dependent subspace, we use a version of Caffarelli's Contraction Principle in Section~\ref{secDependentGaussian}.
Finally, all ingredients are pieced together to prove Theorem~\ref{RBLtheoequa} in Section~\ref{secFinalProof}.

As applications of the understanding the equality case of the Brascamp-Lieb and Barthe's inequalities,
we discuss the Bollobas-Thomason inequality and in its dual version in Section~\ref{secBT0}, and provide the characterization of the equality cases in  Section~\ref{secBT}.

\section{Some applications: Equality in the Bollobas-Thomason inequality and in its dual}
\label{secBT0}

We write $e_1,\ldots,e_n$ to denote an orthonomal basis of $\R^n$.
 For a compact set $K\subset\R^n$ with
${\rm dim}\,{\rm aff}\,K=m$, we write $|K|$ to denote the $m$-dimensional Lebesgue measure of $K$.

The starting point of this section is the classical Loomis-Whitney inequality
\cite{LoW49}.

\begin{theo}[Loomis, Whitney]
\label{Loomis-Whitney}
If $K\subset \R^n$ is compact and affinely spans $\R^n$, then
\begin{equation}
\label{Loomis-Whitney-ineq}
|K|^{n-1}\leq \prod_{i=1}^k|P_{e_i^\bot}K|,
\end{equation}
with equality if and only if
$K=\oplus_{i=1}^nK_i$ where ${\rm aff}K_i$ is a line parallel to $e_i$.
\end{theo}

Meyer \cite{Mey88} provided a dual form of the Loomis-Whitney inequality where equality holds for affine crosspolytopes.

\begin{theo}[Meyer]
\label{dual-Loomis-Whitney}
If $K\subset \R^n$ is compact convex with $o\in{\rm int}K$, then
\begin{equation}
\label{dual-Loomis-Whitney-ineq}
|K|^{n-1}\geq \frac{n!}{n^n}\prod_{i=1}^k|K\cap e_i^\bot|,
\end{equation}
with equality if and only if
$K={\rm conv}\{\pm\lambda_ie_i\}_{i=1}^n$ for $\lambda_i>0$, $i=1,\ldots,n$.
\end{theo}

We note that various Reverse and dual Loomis-Whitney type inequalities are proved by
Campi, Gardner, Gronchi \cite{CGG16}, Brazitikos {\it et al} \cite{BDG17,BGL18},
Alonso-Guti\'errez {\it et al} \cite{ABBC,AB}.

To consider a genarization of the Loomis-Whitney inequality and its dual form,
we set
$[n]:=\{1,\ldots,n\}$, and for a non-empty proper subset $\sigma\subset[n]$, we define
$E_\sigma={\rm lin}\{e_i\}_{i\in\sigma}$. For $s\geq 1$, we say that the not necessarily distinct proper
non-empty subsets $\sigma_1,\ldots,\sigma_k\subset[n]$ form an $s$-uniform cover of $[n]$ if each
$j\in[n]$ is contained in exactly $s$ of $\sigma_1,\ldots,\sigma_k$.

The Bollobas-Thomason inequality \cite{BoT95} reads as follows.

\begin{theo}[Bollobas, Thomason]
\label{Bollobas-Thomason}
If $K\subset \R^n$ is compact and affinely spans $\R^n$, and
$\sigma_1,\ldots,\sigma_k\subset[n]$ form an $s$-uniform cover of $[n]$ for $s\geq 1$, then
\begin{equation}
\label{Bollobas-Thomasson-ineq}
|K|^s\leq \prod_{i=1}^k|P_{E_{\sigma_i}}K|.
\end{equation}
\end{theo}

We note that additional the case when $k=n$, $s=n-1$, and hence when we may assume that $\sigma_i=[n]\backslash e_i$, is the
Loomis-Whitney inequality Therem~\ref{Loomis-Whitney}.

Liakopoulos \cite{Lia19} managed to prove a dual form of the Bollobas-Thomason inequality.
For a finite set $\sigma$, we write $|\sigma|$ to denote its cardinality.

\begin{theo}[Liakopoulos]
\label{Liakopoulos}
If $K\subset \R^n$ is compact convex with $o\in{\rm int}K$, and
$\sigma_1,\ldots,\sigma_k\subset[n]$ form an $s$-uniform cover of $[n]$ for $s\geq 1$, then
\begin{equation}
\label{Liakopoulos-ineq}
|K|^s\geq \frac{\prod_{i=1}^k|\sigma_i|!}{(n!)^s}\cdot \prod_{i=1}^k|K\cap E_{\sigma_i}|.
\end{equation}
\end{theo}

The equality case of
the  Bollobas-Thomason inequality Theorem~\ref{Bollobas-Thomason} based on Valdimarsson \cite{Val08} has been known to the experts but we present this argument in order to have a written account.
Let $s\geq 1$, and let $\sigma_1,\ldots,\sigma_k\subset[n]$ be an $s$-uniform cover
   of $[n]$. We say that
 $\tilde{\sigma}_1,\ldots,\tilde{\sigma}_l\subset[n]$
form a $1$-uniform cover of $[n]$ induced by
the $s$-uniform cover $\sigma_1,\ldots,\sigma_k$ if
$\{\tilde{\sigma}_1,\ldots,\tilde{\sigma}_l\}$ consists of all non-empty distinct subsets of $[n]$
of the form $\cap_{i=1}^k\sigma^{\varepsilon(i)}_i$ where $\varepsilon(i)\in\{0,1\}$ and $\sigma_i^0=\sigma_i$ and  $\sigma_i^1=[n]\setminus\sigma_i$. We observe that $\tilde{\sigma}_1,\ldots,\tilde{\sigma}_l\subset[n]$
actually form a $1$-uniform cover of $[n]$; namely, $\tilde{\sigma}_1,\ldots,\tilde{\sigma}_l$ is a partition of $[n]$.

\begin{theo}[Folklore]
\label{Bollobas-Thomason-eq}
Let $K\subset \R^n$ be compact and affinely span $\R^n$, and let
$\sigma_1,\ldots,\sigma_k\subset[n]$ form an $s$-uniform cover of $[n]$ for $s\geq 1$.
Then equality holds in \eqref{Bollobas-Thomasson-ineq} if and only if
$K=\oplus_{i=1}^l P_{E_{\tilde{\sigma}_i}}K$
where $\tilde{\sigma}_1,\ldots,\tilde{\sigma}_l$ is
the $1$-uniform cover  of $[n]$ induced by
 $\sigma_1,\ldots,\sigma_k$.
\end{theo}

Our main result in this section is the characterization
of the equality case of  the dual Bollobas-Thomason inequality Theorem~\ref{Liakopoulos}
relating it to the Geometric Barthe's inequality.

\begin{theo}
\label{Liakopoulos-eq}
Let $K\subset \R^n$ be compact convex  with $o\in{\rm int}K$, and let
$\sigma_1,\ldots,\sigma_k\subset[n]$ form an $s$-uniform cover of $[n]$ for $s\geq 1$.
Then equality holds in \eqref{Liakopoulos-ineq} if and only if
$K={\rm conv}\{K\cap F_{\tilde{\sigma}_i}\}_{i=1}^l$
where $\tilde{\sigma}_1,\ldots,\tilde{\sigma}_l$ is
the $1$-uniform cover  of $[n]$ induced by
 $\sigma_1,\ldots,\sigma_k$.
\end{theo}


\section{The determinantal inequality and structure theory
for rank one Geometric Brascamp-Lieb data}
\label{secStructure}

We first discuss the basic properties of a set of vectors $u_1,\ldots,u_k\in S^{n-1}$ and constants $c_1,\ldots,c_k>0$
occurring in the Geometric Brascamp-Lieb inequality; namely,
satisfying
\begin{equation}
\label{dim1sum}
\sum_{i=1}^kc_i u_i\otimes u_i=I_n.
\end{equation}
This  section just retells the story of Section~2 of Barthe \cite{Bar98} in the language of
Carlen,  Lieb,  Loss \cite{CLL04} and
 Bennett, Carbery, Christ, Tao \cite{BCCT08}.

\begin{lemma}
\label{dim1properties0}
For  $u_1,\ldots,u_k\in S^{n-1}$ and $c_1,\ldots,c_k>0$ satisfying \eqref{dim1sum},
we have
\begin{description}
\item{(i)} $\sum_{i=1}^k c_i=n$;
\item{(ii)} $\sum_{i=1}^k c_i\langle u_i,x\rangle^2=\|x\|^2$ for all $x\in\R^n$;
\item{(iii)} $c_i\leq 1$ for $i=1,\ldots,k$ with equality if and only if $u_j\in u_i^\bot$ for $j\neq i$;
\item{(iv)} $u_1,\ldots,u_k$ spans $\R^n$, and $k=n$ if and only if $u_1,\ldots,u_n$ is an orthonormal basis of $\R^n$ and $c_1=\ldots=c_n=1$;
\item{(v)} if $L$ is a proper linear subspace of $\R^n$, then
$$
\sum_{u_i\in L}c_i\leq \dim L,
$$
with equality if and only if $\{u_1,\ldots,u_k\}\subset L\cup L^\bot$.
\end{description}
\end{lemma}
\noindent {\bf Remark} If $\sum_{u_i\in L}c_i=\dim L$ in (v), then ${\rm lin}\{u_i:u_i\in L\}=L$
and ${\rm lin}\{u_i:u_i\in L^\bot\}=L^\bot$.\\
\proof Here (i) follows from comparing the traces of the two sides of \eqref{dim1sum}, and (ii) is just an equivalent form of
\eqref{dim1sum}. To prove $c_j\leq 1$ with the characterization of equality, we substitute $x=u_j$ into (ii).

Turning to (iv), $u_1,\ldots,u_k$ spans $\R^n$ by (ii). Next,
let us assume that $u_1,\ldots,u_n\in S^{n-1}$ and $c_1,\ldots,c_n>0$ satisfy \eqref{dim1sum}.
We consider $w_j\in  S^{n-1}$  for $j=1,\ldots,n$ such that  $\langle w_j,u_i\rangle=0$ if $i\neq j$,
and hence (ii) shows that $u_j=\pm w_j$ and $c_j=1$.

For (v), if $u_i\not\in L$, then we consider the unit vector
$$
\tilde{u}_i=\frac{P_{L^\bot}u_i}{\|P_{L^\bot}u_i\|}\in L^\bot.
$$
We deduce that if $x\in L^\bot$, then
$$
\|x\|^2=\sum_{i=1}^k c_i \langle u_i,x\rangle^2=\sum_{u_i\not\in L}c_i\langle P_{L^\bot}u_i,x\rangle^2
=\sum_{u_i\not\in L}c_i\|P_{L^\bot}u_i\|^2\langle \tilde{u}_i,x\rangle^2.
$$
It follows from (i) and (ii)  applied to $\{\tilde{u}_i:\,u_i\not\in L\}$ in $L^\bot$  that
$$
\dim L^\bot=\sum_{u_i\not\in L}c_i\left\|P_{L^\bot}u_i\right\|^2\leq \sum_{u_i\not\in L}c_i.
$$
In turn, we conclude the inequality in (v) by (i). Equality holds in (v) if and only if
$\|P_{L^\bot}u_i\|=1$ whenever $u_i\not\in L$; therefore, $u_1,\ldots,u_k\subset L\cup L^\bot$.
\proofbox

Let $u_1,\ldots,u_k\in S^{n-1}$ and $c_1,\ldots,c_k>0$ satisfy \eqref{dim1sum}. Following  Bennett, Carbery, Christ, Tao \cite{BCCT08}, we say that a non-zero linear subspace $V$ is a
critical subspace with respect to $u_1,\ldots,u_k$ and  $c_1,\ldots,c_k$ if
$$
\sum_{u_i\in V}c_i=\dim V.
$$
In particular, $\R^n$ is a critical subspace according to  Lemma~\ref{dim1properties0}.
We say that a non-empty subset \ $\mathcal{U}\subset \{u_1,\ldots,u_k\}$ is indecomposable
if ${\rm lin}\,\mathcal{U}$ is an indecomposable critical subspace.

In order to understand the equality case of the rank one Brascamp-Lieb inequality, Barthe \cite{Bar98}
 indicated an equivalence relation on $\{u_1,\ldots,u_k\}$.
We say that a subset $\mathcal{D}\subset\{u_1,\ldots,u_k\}$ is minimally dependent if   $\mathcal{D}$ is dependent and no proper subset of $\mathcal{D}$ is dependent. The following is folklore in matroid theory, was known most probably
already to Tutte (see for example Theorem 7.3.6 in Recski \cite{Rec89}). For the convenience of the reader, we provide an argument.

\begin{lemma}
\label{bowtie}
Given non-zero $v_1,\ldots,v_k$ spanning $\R^n$, $n\geq 1$, we write $v_i\bowtie v_j$
if either $v_i=v_j$, or there exists a minimal dependent set $\mathcal{D}\subset\{v_1,\ldots,v_k\}$ satisfying
$v_i,v_j\in \mathcal{D}$.
\begin{description}
\item{(i)} $v_i\bowtie v_j$ if and only if there exists a subset $\mathcal{U}\subset\{v_1,\ldots,v_k\}$ of cardinality
$n-1$ such that both $\{v_i\}\cup\cal{U}$ and $\{v_j\}\cup\cal{U}$ are independent;
\item{(ii)} $\bowtie$ is an equivalence relation
on $\{v_1,\ldots,v_k\}$;
\item{(iii)} if $V_1,\ldots,V_m$ are the linear hulls of the equivalence classes with respect to $\bowtie$, then
they span $\R^n$ and $V_i\cap V_j=\{0\}$ for $i\neq j$.
\end{description}
\end{lemma}
\proof We prove the lemma  by induction on $n\geq 1$ where the case $n=1$ readily holds. Therefore, we assume that
 $n\geq 2$.

We may readily assume that
\begin{equation}
\label{onedim-unique}
\{v_1,\ldots,v_k\}\cap {\rm lin}\,\{v_i\}=\{v_i\}\mbox{ \ \ for $i=1,\ldots,k$.}
\end{equation}

For (i), if $\mathcal{D}$ is a minimal dependent set with
$v_i,v_j\in \mathcal{D}$, then adding some $\mathcal{V}\subset\{v_1,\ldots,v_k\}$ to $\mathcal{D}\backslash\{v_i\}$,
we obtain a basis of $\R^n$, and we may choose $\mathcal{U}=\mathcal{V}\cup(\mathcal{D}\backslash\{v_i,v_j\})$.
On the other hand, if the suitable $\cal{U}$ of cardinality $n-1$ exists such that both $\{v_i\}\cup\cal{U}$ and
$\{v_j\}\cup\cal{U}$ are independent, then any dependent subset of $\mathcal{U}\cup\{v_i,v_j\}$ contains $v_i$ and $v_j$.

For (ii) and (iii), we call a non-zero linear subspace $W\subset \R^n$ unsplittable with respect to $\{v_1,\ldots,v_k\}$ if $W$ is spanned by
$W\cap\{v_1,\ldots,v_k\}$, but
there exist no non-zero complementary linear subspaces $A,B\subset W$ with
$\{v_1,\ldots,v_k\}\cap W\subset A\cup B$. Readily, there exist pairwise complementary unsplittable linear subspaces
$W_1,\ldots,W_m\subset \R^n$
such that $\{v_1,\ldots,v_k\}\subset W_1\cup \ldots \cup W_m$.

On the one hand, if $v_i\in W_\alpha$ and $v_j\in W_\beta$ for $\alpha\neq \beta$, then trivially $v_i\not\bowtie v_j$. Therefore all we need to prove is that if $v_i,v_j\in W_\alpha$, then $v_i\bowtie v_j$.
By the induction on $n$, we may assume that $m=1$ and $W_\alpha=\R^n$. We may also assume that $i=1$ and $j=2$.

The final part of argument is indirect; therefore, we suppose that
\begin{equation}
\label{nobowtie}
v_1\not\bowtie v_2,
\end{equation}
and seek a contradiction.

\eqref{nobowtie} implies that $v_1$ and $v_2$ are independent, and hence
$v_1\not\bowtie v_2$ and \eqref{onedim-unique} yield that
$L={\rm lin}\{v_1,v_2\}$ satisfies
\begin{equation}
\label{v1v2-only}
\{v_1,\ldots,v_k\}\cap L=\{v_1,v_2\}.
\end{equation}
 Now $\R^n$ is unsplittable, thus $n\geq 3$.

Since $v_1,\ldots,v_k$ span $\R^n$, we may assume that $v_1,\ldots,v_n$ form a basis of $\R^n$.
Let $L_0={\rm lin}\{v_3,\ldots,v_n\}$,  and $L_t={\rm lin}\{v_t,L_0\}$ for
$t=1,2$. We may also assume that $v_1,\ldots,v_n$ is an orthonormal basis.

For any $l>n$, (i) and $v_1\not\bowtie v_2$ yield that
\begin{equation}
\label{eitherL1L2}
\mbox{either $v_l\in L_1$, or $v_l\in L_2$}.
\end{equation}
Since $\R^n$ is unsplittable, there exist $p,q>n$ such that
\begin{equation}
\label{L1vpL2vq}
\mbox{$v_p\in L_1\backslash L_0$ and $v_q\in L_2\backslash L_0$}.
\end{equation}

For any $w\not\in L$, we write
$$
{\rm supp}\,w=\{v_l:\,l\in\{3,\ldots,n\}\;\&\;\langle w,v_l\rangle\neq 0\};
$$
namely, the basis vectors where the corrresponding coordinate of $w|L-0$ is non-zero.\\

\noindent{\bf Case 1} {\it There exist $v_p\in L_1\backslash L_0$ and $v_q\in L_2\backslash L_0$, $p,q>n$, such that
$({\rm supp}\,v_p)\cap({\rm supp}\,v_q)\neq\emptyset$}

Let $v_s\in ({\rm supp}\,v_p)\cap({\rm supp}\,v_q)$.
Now the $n+1$ element set
$$
\{v_1,v_p,v_2,v_q\}\cup\{v_l:\,l\in\{3,\ldots,n\}\backslash \{s\}\}
$$
is dependent,
and considering the $1^{st}$, $2^{nd}$ and $s^{th}$ coordinates
show that both $v_1$ and $v_2$ lie in any dependent subset. This fact contradicts \eqref{nobowtie}.\\

 \noindent{\bf Case 2} {\it $({\rm supp}\,v_p)\cap({\rm supp}\,v_q)=\emptyset$
for any $v_p\in L_1\backslash L_0$ and $v_q\in L_2\backslash L_0$ with $p,q>n$}

Let ${\cal U}_t=\cup\{{\rm supp}\,v_p:\,p>n\;\&\;v_p\in L_t\backslash L_0\}$ for $t=1,2$. It follows that
${\cal U}_1\cap {\cal U}_2=\emptyset$, thus $n\geq 4$. For any partition ${\cal U}'_1\cup {\cal U}'_2=\{v_3,\ldots,v_n\}$
(and hence ${\cal U}'_1\cap {\cal U}'_2=\emptyset$) such that ${\cal U}_1\subset {\cal U}'_1$
and ${\cal U}_2\subset {\cal U}'_2$, there exists some $v_l\in L_0$ that is contained
neither in ${\rm lin}\,({\cal U}'_1\cup\{v_1\})$ nor in ${\rm lin}\,({\cal U}'_2\cup\{v_2\})$ because $\R^n$ is unslittable.
In turn we deduce that we may reindex the vectors $v_3,\ldots,v_n$ on the one hand, and the vectors
$v_{n+1},\ldots,v_k$ on the other hand to ensure the following properties:
\begin{itemize}
\item $v_{n+1}\in L_1\backslash L_0$ and $v_{n+2}\in L_2\backslash L_0$;
\item there exist $\alpha\in\{3,\ldots,n-1\}$ and $\beta\in\{n+3,\ldots,k\}$ such that
${\rm supp}\,v_l\subset\{v_\alpha,\ldots,v_n\}$ for $l\in\{n+1,\ldots,\beta\}$, and $v_l\in L_0$ if
$n+3\leq l\leq\beta$;
\item for any partion ${\cal W}_1\cup{\cal W}_2=\{v_\alpha,\ldots,v_n\}$ into non-empty sets,
there exist $l\in\{n+1,\ldots,\beta\}$ such that ${\rm supp}\,v_l$ intersects both
${\cal W}_1$ and ${\cal W}_2$.
\end{itemize}
We observe that $\widetilde{L}_0={\rm lin}\{v_\alpha,\ldots,v_n\}$ is unsplittable with respect to
$$
\{v_\alpha,\ldots,v_n,v_{n+1}|L_0,v_{n+2}|L_0,v_{n+3},\ldots,v_\beta\}.
$$
Therefore, this last set contains a minimal dependent subset $\widetilde{{\cal D}}$ with
$v_{n+1}|L_0,v_{n+2}|L_0\in \widetilde{{\cal D}}$ by induction; namely, the elements of $\widetilde{{\cal D}}$
different from $v_{n+1}|L_0,v_{n+2}|L_0$ are vectors of the form $v_l$ that lie in $L_0$. We conclude that
$$
{\cal D}=\{v_1,v_2,v_{n+1},v_{n+2}\}\cup
\left(\widetilde{{\cal D}}\backslash\{v_{n+1}|L_0,v_{n+2}|L_0\}\right)
$$
is a minimal dependent set, contradicting \eqref{nobowtie}, and proving Lemma~\ref{bowtie}.
\proofbox

\begin{lemma}
\label{dim1properties1}
For  $u_1,\ldots,u_k\in S^{n-1}$ and $c_1,\ldots,c_k>0$ satisfying \eqref{dim1sum},
we have
\begin{description}
\item{(i)} a proper linear subspace $V\subset \R^n$ is critical if and only if $\{u_1,\ldots,u_k\}\subset V\cup V^\bot$;
\item{(ii)} if $V,W$ are proper critical subspaces with $V\cap W\neq\{0\}$, then $V^\bot$, $V\cap W$ and $V+W$
are critical subspaces;
\item{(iii)} the equivalence classes with respect to the relation $\bowtie$ in Lemma~\ref{bowtie} are the indecomposable subsets
of $\{u_1,\ldots,u_k\}$;
\item{(iv)} the proper indecomposable critical subspaces are pairwise orthogonal, and any critical subspace is the sum of some
 indecomposable critical subspaces.
 \end{description}
\end{lemma}
\proof (i) directly follows from Lemma~\ref{dim1properties0} (v), and in turn (i) yields (ii).

We prove (iii) and and  first half of (iv) simultatinuously. Let $V_1,\ldots,V_m$ be
the linear hulls of the equivalence classes of $u_1,\ldots,u_k$ with respect to the $\bowtie$ of Lemma~\ref{bowtie}.
We deduce from Lemma~\ref{dim1properties0} (v) that each $V_i$ is a critical subspace, and
if $i\neq j$, then $V_i$ and $V_j$ are orthogonal.

Next let $\mathcal{U}\subset\{u_1,\ldots,u_k\}$ be an indecomposable set, and let $V={\rm lin}\,\mathcal{U}$.
We write $I\subset \{1,\ldots,m\}$ to denote the set of indices $i$ such that $V_i\cap \mathcal{U}\neq \emptyset$.
Since $V$ is a critical subspace, we deduce from Lemma~\ref{dim1properties0} (v) that $V_i\cap V$ is a critical subspace
for $i\in I$, as well; therefore, $I$ consists of a unique index $p$ as $\mathcal{U}$ is indecomposable.
In particular, $V=V_p$.

It follows from Lemma~\ref{dim1properties0} (v) that $\{u_1,\ldots,u_k\}\subset V\cup V^\bot$; therefore, there exists
no minimally dependent subset of $\{u_1,\ldots,u_k\}$ intersecting both $\mathcal{U}$ and its complement. We conclude that
$V=V_p$.

Finally, the second half of (iv)  follows from (i) and (ii).
\proofbox

The following is the main result of this section, where the inequality is proved by
Ball \cite{Bal89,Bal91}, and the equality case is clarified by Barthe \cite{Bar98}.

\begin{prop}[Ball-Barthe Lemma]
\label{GBLconstantdim1}
For  $u_1,\ldots,u_k\in S^{n-1}$ and $c_1,\ldots,c_k>0$ satisfying \eqref{dim1sum},
if $t_i>0$  for $i=1,\ldots,k$, then
\begin{equation}
\label{BallBarthe}
\det\left( \sum_{i=1}^kc_it_iu_i\otimes u_i\right)\geq \prod_{i=1}^k t_i^{c_i}.
\end{equation}
Equality holds in (\ref{BallBarthe})
if and only if $t_i=t_j$ for any $u_i$ and $u_j$ lying in the same indecomposable subset
of $\{u_1,\ldots,u_k\}$.
\end{prop}
\proof To simplify expressions, let $v_i=\sqrt{c_i}u_i$ for $i=1,\ldots,k$.

In this argument, $I$ always denotes some subset of $\{1,\ldots,k\}$ of cardinality $n$.
For $I=\{i_1,\ldots,i_n\}$, we define
$$
d_I:=  \det[v_{i_1},\ldots,v_{i_n}]^2\qquad \text{and}\qquad
t_I := t_{i_1}\cdots t_{i_n}.
$$
For the $n\times k$ matrices $M=[v_1,\ldots,v_k]$ and
$\widetilde{M}=[\sqrt{t_1}\,v_1,\ldots,\sqrt{t_k}\,v_k]$, we have
\begin{equation}
M M^T=I_n\mbox{ \ and \ }\widetilde{M} \widetilde{M}^T=\sum_{i=1}^kt_iv_i\otimes v_i.
\end{equation}
It follows from
 the Cauchy-Binet formula that
$$
\sum_I d_I= 1\qquad \text{and}\qquad
\det \left(\sum_{i=1}^kt_iv_i\otimes v_i\right) =\sum_I t_Id_I,
$$
where the summations extend over all sets $I\subset\{1,\ldots,k\}$ of cardinality $n$.
It follows that the discrete measure $\mu$ on the
$n$ element subsets of $\{1,\ldots,k\}$ defined by $\mu(\{I\})=d_I$
is a probability measure. We deduce from inequality between the arithmetic and geometric mean that
\begin{equation}
\label{agmean}
\det\left( \sum_{i=1}^kt_iv_i\otimes v_i\right) =\sum_I t_Id_I\geq
\prod_It_I^{d_I}.
\end{equation}

The factor $t_i$ occurs in $\prod_It_I^{d_I}$ exactly
$\sum_{I,\,i\in I}d_I$ times. Moreover, the Cauchy-Binet formula applied
to the vectors $v_1,\ldots,v_{i-1},v_{i+1},\ldots,v_k$ implies
\begin{align*}
\sum_{I,\,i\in I}d_I&=\sum_Id_I-\sum_{I,\,i\not\in I}d_I=
1-\det\left(\sum_{j\neq i}v_j\otimes v_j\right)\\
&=1-\det\left({\rm Id}_n-v_i\otimes v_i\right)=\langle v_i,v_i\rangle=c_i.
\end{align*}
Substituting this into (\ref{agmean}) yields (\ref{BallBarthe}).

We now assume that equality holds in (\ref{BallBarthe}).
Since equality holds in (\ref{agmean}) when applying
arithmetic and geometric mean,  all the $t_I$  are the same
for any subset $I$ of $\{1,\ldots,k\}$ of cardinality $n$ with $d_I\neq 0$.
It follows that $t_i=t_j$ whenever $u_i\bowtie u_j$, and in turn we deduce
that $t_i=t_j$ whenever $u_i$ and $u_j$ lie in the same indecomposable set
by Lemma~\ref{dim1properties1} (i).

On the other hand, Lemma~\ref{dim1properties1} (ii) yields that
if $t_i=t_j$ whenever $u_i$ and $u_j$ lie in the same indecomposable set,
then equality holds in (\ref{BallBarthe}).
\proofbox

Combining Lemma~\ref{dim1properties1} and Proposition~\ref{GBLconstantdim1} leads to the following:

\begin{coro}
\label{GBLconstantdim1cor}
For  $u_i\in S^{n-1}$ and $c_i,t_i>0$, $i=1,\ldots,k$ satisfying \eqref{dim1sum},
equality holds in (\ref{BallBarthe})
if and only if there exist pairwise orthogonal linear subspaces $V_1,\ldots,V_m$, $m\geq 1$, such that
$\{u_1,\ldots,u_k\}\subset V_1\cup\ldots\cup V_m$ and
$t_i=t_j$ whenever $u_i$ and $u_j$ lie in the same $V_p$ for some $p\in\{1,\ldots,m\}$.
\end{coro}

\section{Structure theory of a Brascamp-Lieb data and
the determinantal inequality corresponding to the higher rank case}
\label{secStuctureHigh}

We build a structural theory for a Brascamp-Lieb data based on results proved or indicated in Barthe \cite{Bar98},
Bennett, Carbery, Christ, Tao \cite{BCCT08} and
 Valdimarsson \cite{Val08}.

For non-zero linear subspaces $E_1,\ldots,E_k$ of $\R^n$ and $c_1,\ldots,c_k>0$  satisfying
the Geometric Brascamp-Lieb condition
\begin{equation}
\label{highdimcond}
\sum_{i=1}^kc_iP_{E_i}=I_n,
\end{equation}
we connect \eqref{highdimcond} to \eqref{dim1sum}.
For $i=1,\ldots,k$, let ${\rm dim}\,E_i=n_i$ and let $u_1^{(i)},\ldots,u_{n_i}^{(i)}$
be any orthonormal basis  of $E_i$. In addition,
for $i=1,\ldots,k$,
we consider the $n \times n_i$ matrix
$M_i=\sqrt{c_i}[u_1^{(i)},\ldots,u_{n_i}^{(i)}]$. We deduce that
\begin{eqnarray}
\label{Eiui}
c_iP_{E_i}&=&M_iM_i^T=\sum_{j=1}^{n_i}c_iu_j^{(i)}\otimes u_j^{(i)}
\mbox{ \ for $i=1,\ldots,k$};\\
\label{high-1}
I_n&=&\sum_{i=1}^kc_iP_{E_i}=
\sum_{i=1}^k\sum_{j=1}^{n_i}c_i u_j^{(i)}\otimes u_j^{(i)}
=\sum_{i=1}^k\sum_{j=1}^{n_i}c_j^{(i)} u_j^{(i)}\otimes u_j^{(i)}
\end{eqnarray}
and hence $u_j^{(i)}\in S^{n-1}$ and $c_j^{(i)}=c_i>0$ for $i=1,\ldots,k$ and
$j=1,\ldots,n_i$ form a
Geometric Brascamp-Lieb data like in \eqref{dim1sum}.

\begin{lemma}
\label{highdimcritical}
For linear subspaces $E_1,\ldots,E_k$ of $\R^n$ and $c_1,\ldots,c_k>0$ satisfying
\eqref{highdimcond},
\begin{description}

\item{(i)} if $x\in \R^n$, then $\sum_{i=1}^kc_i\|P_{E_i}x\|^2=\|x\|^2$;

\item{(ii)} if $V\subset \R^n$ is a proper linear subspace, then
\begin{equation}
\label{highdimcondeq}
\sum_{i=1}^kc_i\dim(E_i\cap V)\leq\dim V
\end{equation}
where equality holds if and only if $E_i=(E_i\cap V)+ (E_i\cap V^\bot)$ for $i=1,\ldots,k$; or equivalently,
when $V=(E_i\cap V)+ (E_i^\bot\cap V)$
for $i=1,\ldots,k$.
\end{description}
\end{lemma}
\proof For $i=1,\ldots,k$, let ${\rm dim}\,E_i=n_i$ and let $u_1^{(i)},\ldots,u_{n_i}^{(i)}$
be any orthonormal basis  of $E_i$ such that if $V\cap E_i\neq\{0\}$, then
$u_1^{(i)},\ldots,u_{m_i}^{(i)}$
is any orthonormal basis  of $V\cap E_i$ where $m_i\leq n_i$.

For any $x\in\R^n$ and $i=1,\ldots,k$, we have $\|P_{E_i}x\|^2=\sum_{j=1}^{n_i}\langle u_j^{(i)},x\rangle^2$,
thus Lemma~\ref{dim1properties0} (ii) yields (i).

Concerning (ii), Lemma~\ref{dim1properties0} (v) yields \eqref{highdimcondeq}.
On the other hand, if equality holds
in \eqref{highdimcondeq}, then $V$ is a critical subspace
for the rank one Geometric Brascamp-Lieb data
$u_j^{(i)}\in S^{n-1}$ and $c_j^{(i)}=c_i>0$ for $i=1,\ldots,k$ and
$j=1,\ldots,n_i$ satisfying \eqref{high-1}.
Thus Lemma~\ref{highdimcritical} (ii) follows from Lemma~\ref{dim1properties0} (v).
\proofbox

We say that a non-zero linear subspace $V$ is a
critical subspace with respect to  the proper linear subspaces  $E_1,\ldots,E_k$ of $\R^n$ and $c_1,\ldots,c_k>0$ satisfying
\eqref{highdimcond} if
$$
\sum_{i=1}^kc_i\dim(E_i\cap V)=\dim V.
$$
In particular, $\R^n$ is a critical subspace by calculating traces of both sides of \eqref{highdimcond}.
For a proper linear subspace $V\subset \R^n$,
Lemma~\ref{highdimcritical} yields that $V$ is critical if and only if $V^\bot$ is critical, which is turn equivalent saying that
\begin{equation}
\label{highdimcritical0}
\mbox{$E_i=(E_i\cap V)+ (E_i\cap V^\bot)$
for $i=1,\ldots,k$;}
\end{equation}
or in other words,
\begin{equation}
\label{highdimcriticalV}
\mbox{$V=(E_i\cap V)+ (E_i^\bot\cap V)$
for $i=1,\ldots,k$.}
\end{equation}
We observe that \eqref{highdimcritical0} has the following consequence: If $V_1$ and $V_2$ are orthogonal critical subspaces,
 then
\begin{equation}
\label{highdimcriticalV1V2}
E_i\cap(V_1+V_2)=(E_i\cap V_1)+ (E_i\cap V_2)
\mbox{ \ for $i=1,\ldots,k$.}
\end{equation}

We recall that a critical subspace $V$ is indecomposable if
$V$ has no proper critical linear subspace.

\begin{lemma}
\label{highdimcritical-cap}
If $E_1,\ldots,E_k$ are linear subspaces of $\R^n$ and $c_1,\ldots,c_k>0$ satisfying
\eqref{highdimcond}, and
 $V,W$ are proper critical subspaces, then $V^\bot$ and $V+W$
are critical subspaces, and even $V\cap W$ is critical provided that $V\cap W\neq\{0\}$.
\end{lemma}
\proof We may assume that ${\rm dim}\,E_i\geq 1$ for $i=1,\ldots,k$.

The fact that $V^\bot$ is also critical follows directly from \eqref{highdimcritical0}.

Concerning $V\cap W$ when $V\cap W\neq\{0\}$, we need to prove that if $i=1,\ldots,k$, then
\begin{equation}
\label{VcapW}
(V\cap W)\cap E_i+(V\cap W)^\bot\cap E_i=E_i.
\end{equation}
For a linear subspace $L\subset E_i$, we write $L^{\bot_i}=L^\bot\cap E_i$ to denote the orthogonal
complement within $E_i$. We observe that as $V$ and $W$ are critical subspaces, we have
$(V\cap E_i)^{\bot_i}= V^\bot \cap E_i$ and $(W\cap E_i)^{\bot_i}= W^\bot \cap E_i$.
It follows from the identity
 $(V\cap W)^\bot=V^\bot+W^\bot$ that
\begin{eqnarray*}
E_i&\supset&(V\cap W)\cap E_i+(V\cap W)^\bot\cap E_i=(V\cap E_i)\cap (W\cap E_i)+(V^\bot+ W^\bot)\cap E_i\\
&\supset&(V\cap E_i)\cap (W\cap E_i)+(V^\bot \cap E_i)+(W^\bot\cap E_i)\\
&=&
(V\cap E_i)\cap (W\cap E_i)+(V\cap E_i)^{\bot_i}+(W\cap E_i)^{\bot_i}\\
&=&
(V\cap E_i)\cap (W\cap E_i)+[(V\cap E_i)\cap (W\cap E_i)]^{\bot_i}=E_i,
\end{eqnarray*}
yielding \eqref{VcapW}.

Finally, $V+W$ is also critical as $V+W=(V^\bot\cap W^\bot)^\bot$.
\proofbox

We deduce from Lemma~\ref{highdimcritical-cap} that any critical subspace can be decomposed into indecomposable ones.

\begin{coro}
\label{critical-decompose}
If $E_1,\ldots,E_k$ are proper linear subspaces of $\R^n$ and $c_1,\ldots,c_k>0$ satisfy
\eqref{highdimcond}, and  $W$ is a critical subspace or $W=\R^n$, then there exist pairwise orthogonal indecomposable critical subspaces
$V_1,\ldots,V_m$, $m\geq 1$, such that $W=V_1+\ldots+V_m$ (possibly $m=1$ and $W=V_1$).
\end{coro}

We note that the decomposition
 of $\R^n$ into indecomposable critical subspaces is not unique in general
 for a Geometric Brascamp-Lieb data.  Valdimarsson \cite{Val08} provides some examples, and in addition, we provide an example where we have a continuous family of indecomposable critical subspaces.

\begin{example}[Continuous family of indecomposable critical subspaces]
\label{manycritical}
In $\R^4$, let us consider the following six unit vectors:
$u_1(1,0,0,0)$,
$u_2(\frac12,\frac{\sqrt{3}}2,0,0)$,
$u_3(\frac{-1}2,\frac{\sqrt{3}}2,0,0)$  ,
$v_1(0,0 ,1,0)$,
$v_2(0,0 ,\frac12,\frac{\sqrt{3}}2)$,
$v_3(0,0,\frac{-1}2,\frac{\sqrt{3}}2)$,
which satisfy $u_2=u_1+u_3$  and $v_2=v_1+v_3$.

For any $x\in \R^4$, we have
$$
\|x\|^2=\sum_{i=1}^3 \frac23 \cdot (\langle x,u_i\rangle^2
+\langle x,v_i\rangle^2 )
$$

  Therefore, we define the Geometric Brascamp-Lieb Data $E_i={\rm lin}\{u_i,v_i\}$ and $c_i=\frac23$ for $i=1,2,3$
satisfying (\ref{highdimcond0}). In this case, $F_{\rm dep}=\R^4$.

For any angle $t\in\R$, we have a two-dimensional indecomposable critical subspace
$$
V_t={\rm lin}\{(\cos t)u_1+(\sin t)v_1,  (\cos t)u_2+(\sin t)v_2,
 (\cos t)u_3+(\sin t)v_3 \}.
$$
\end{example}

Next we prove the crucial determinantal inequality. Its proof is kindly provided by Franck Barthe.

\begin{prop}[Barthe]
\label{GBLconstanthighdim}
For linear subspaces $E_1,\ldots,E_k$ of $\R^n$, $n\geq 1$ and $c_1,\ldots,c_k>0$ satisfying
\eqref{highdimcond},
if $A_i:E_i\to E_i$ is a positive definite linear transformation for $i=1,\ldots,k$, then
\begin{equation}
\label{BallBarthehigh}
\det\left( \sum_{i=1}^kc_iA_iP_{E_i}\right)\geq \prod_{i=1}^k(\det A_i)^{c_i}.
\end{equation}
Equality holds in \eqref{BallBarthehigh}
if and only if
there exist linear subspaces $V_1,\ldots,V_m$ where $V_1=\R^n$ if $m=1$ and
$V_1,\ldots,V_m$ are pairwise orthogonal indecomposable critical subspaces spanning $\R^n$
if $m\geq 2$, and a positive definite $n\times n$ matrix $\Phi$ such that
$V_1,\ldots,V_m$ are eigenspaces of $\Phi$ and $\Phi|_{E_i}=A_i$ for $i=1,\ldots,k$.
In addition, $\Phi=\sum_{i=1}^k c_iA_iP_{E_i}$ in the case of equality.
\end{prop}
\proof We may assume that ${\rm dim}\,E_i\geq 1$ for $i=1,\ldots,k$.

For $i=1,\ldots,k$, let ${\rm dim}\,E_i=n_i$, let $u_1^{(i)},\ldots,u_{n_i}^{(i)}$
be an orthonormal basis  of $E_i$ consisting of eigenvectors of $A_i$,
and let $\lambda^{(i)}_j>0$ be the eigenvalue of $A_i$ corresponding to  $u_j^{(i)}$.
In particular $\det A_i=\prod_{j=1}^{n_i}\lambda^{(i)}_j$ for $i=1,\ldots,k$.
 In addition,
for $i=1,\ldots,k$,
we set
$M_i=\sqrt{c_i}[u_1^{(i)},\ldots,u_{n_i}^{(i)}]$ and $B_i$ to be the positive definite
transformation with $A_i=B_iB_i$, and hence
$$
c_iA_iP_{E_i}=(M_iB_i)(M_iB_i)^T=\sum_{j=1}^{n_i}c_i\lambda^{(i)}_ju_j^{(i)}\otimes u_j^{(i)}.
$$

We deduce  from Lemma~\ref{GBLconstantdim1} and \eqref{high-1} that
\begin{eqnarray}
\nonumber
\det\left( \sum_{i=1}^kc_iA_iP_{E_i}\right)&=&
\det\left( \sum_{i=1}^k\sum_{j=1}^{n_i}c_i\lambda^{(i)}_ju_j^{(i)}\otimes u_j^{(i)}\right)\\
\label{GBLconstantdim10}
&\geq &
\prod_{i=1}^k\left(\prod_{j=1}^{n_i} \lambda^{(i)}_j\right)^{c_i}
=\prod_{i=1}^k(\det A_i)^{c_i}.
\end{eqnarray}

If we have equality in \eqref{BallBarthehigh}, and hence also in
\eqref{GBLconstantdim10}, then Corollary~\ref{GBLconstantdim1cor}
implies that there exist pairwise orthogonal critical subspaces $V_1,\ldots,V_m$, $m\geq 1$ spanning $\R^n$
and $\lambda_1,\ldots,\lambda_m>0$ (where $V_1=\R^n$ if $m=1$) such that
if $E_i\cap V_j\neq\{0\}$, then $E_i\cap V_j$ is an eigenspace of $A_i$ with eigenvalue $\lambda_j$.
We conclude from \eqref{highdimcritical0} that each $V_j$ is a critical subspace,
and from Corollary~\ref{critical-decompose} that each $V_j$ can be assumed to be indecomposable.
Finally, \eqref{highdimcriticalV1V2} yields that each $E_i$ is spanned by
the subspaces $E_i\cap V_j$ for $j=1,\ldots,m$.

To show that each $V_j$ is an eigenspace for the positive definite linear transform
$\sum_{i=1}^k c_iA_iP_{E_i}$ of $\R^n$ with eigenvalue $\lambda_j$, we observe that
$$
A_iP_{E_i}x=\lambda_jP_{E_i}x
$$
for any $i=1,\ldots,k$ and $x\in V_j$. It follows that if $x\in V_j$, then
$$
\sum_{i=1}^k c_iA_iP_{E_i}x=\lambda_j\sum_{i=1}^k c_iP_{E_i}x=\lambda_jx,
$$
proving that we can choose $\Phi=\sum_{i=1}^k c_iA_iP_{E_i}$.

On the other hand, let us assume that there exists
a positive definite $n\times n$ matrix $\Theta$ whose
eigenspaces $W_1,\ldots,W_l$ are critical subspaces (or $l=1$ and $W_1=\R^n$) and $\Theta|_{E_i}=A_i$ for $i=1,\ldots,k$.
In this case, for any $i=1,\ldots,k$, we may choose the orthonormal basis $u_1^{(i)},\ldots,u_{n_i}^{(i)}$ of $E_i$
in a way such that $u_1^{(i)},\ldots,u_{n_i}^{(i)}\subset W_1\cup\ldots\cup W_l$, and hence
Corollary~\ref{GBLconstantdim1cor} yields that equality holds in \eqref{BallBarthehigh}.
\proofbox

\noindent{\bf Remark } While Proposition~\ref{GBLconstanthighdim} has a crucial role in proving both the Brascamp-Lieb inequality \eqref{BL} and Barthe's inequality \eqref{RBL} and their equality cases,
Proposition~\ref{GBLconstanthighdim} can be actually derived from say \eqref{BL}. In the Brascamp-Lieb inequality,
choose $f_i(z)=e^{-\pi\langle A_i z,z \rangle}$ for $z\in E_i$ and $i=1,\ldots,k$, and hence
$\int_{E_i}f_i=\left(\det A_i\right)^{\frac{-1}2}$. On the other hand, if $x\in\R^n$, then
$$
\prod_{i=1}^k f_i\left(P_{E_i}x\right)^{c_i}=e^{-\pi\sum_{i=1}^kc_i\langle A_i P_{E_i}x,P_{E_i}x \rangle}
=e^{-\pi\sum_{i=1}^kc_i\langle A_i P_{E_i}x,x \rangle}
=e^{-\pi\langle\sum_{i=1}^kc_i A_i P_{E_i}x,x \rangle};
$$
therefore, the Brascamp-Lieb inequality \eqref{BL} yields
$$
\left(\det \sum_{i=1}^kc_i A_i P_{E_i}\right)^{\frac{-1}2}\leq \prod_{i=1}^k\left(\det A_i\right)^{\frac{-c_i}2}.
$$
In addition, the equality conditions in Proposition~\ref{GBLconstanthighdim} can be derived from Valdimarsson's Theorem~\ref{BLtheoequa}.\\

Let us show why indecomposability of the critical subspaces in Proposition~\ref{GBLconstanthighdim} is useful.

\begin{lemma}
\label{indec-dependent-critical}
Let the linear subspaces $E_1,\ldots,E_k$ of $\R^n$ and $c_1,\ldots,c_k>0$ satisfy
\eqref{highdimcond}, let $F_{\rm dep}\neq \R^n$, and let $F_1,\ldots,F_l$ be the independent subspaces, $l\geq 1$.
If $V$ is an indecomposable critical subspace, then either $V\subset F_{\rm dep}$, or there exists an independent
subspace $F_j$, $j\in\{1,\ldots,l\}$ such that $V\subset F_j$.
\end{lemma}
\proof It is equivalent to prove that if $V$ is an indecomposable critical subspace and $j\in\{1,\ldots,l\}$, then
\begin{equation}
\label{FjVnonintersect}
V\not \subset F_j\mbox{ \ implies \ } F_j\subset V^\bot.
\end{equation}
We deduce that $V\cap F_j=\{0\}$ from the facts that $V$ is indecomposable
and $F_j$ is a critical subspace, thus $F_j\cap V$ is a critical subspace or
$\{0\}$.
There exists a partion $M\cup N=\{1,\ldots,k\}$ with $M\cap N=\emptyset$ such that
$$
F_j=\left(\cap_{i\in M}E_i\right)\cap \left(\cap_{i\in N}E_i^\bot\right).
$$

Let $y\in F_j$.
Since $V$ is a critical subspace, we conclude that $P_Vy\in E_i$ for $i\in M$ and $P_Vy\in E_i^\bot$ for $i\in N$,
and hence $P_Vy\in V\cap \left(\cap_{i\in M}E_i\right)\cap \left(\cap_{i\in N}E_i^\bot\right)=\{0\}$.
Therefore, $y\in V^\bot$.
\proofbox

\section{Typical Gaussian extremizers for some Geometric Brascamp-Lieb data}
\label{secGaussians}

This section continues to build on work done in Barthe \cite{Bar98},
Bennett, Carbery, Christ, Tao \cite{BCCT08} and
 Valdimarsson \cite{Val08}.

For linear subspaces $E_1,\ldots,E_k$ of $\R^n$ and $c_1,\ldots,c_k>0$ satisfying
\eqref{highdimcond}, we deduce from Lemma~\ref{highdimcritical} (i) and \eqref{highdimcriticalV} that
if $V$ is a critical subspace, then writing $P^{(V)}_{E_i\cap V}$ to denote the restriction of
$P_{E_i\cap V}$ onto $V$, we have
\begin{equation}
\label{highdimcriticalsum}
\sum_{E_i\cap V\neq\{0\}}c_iP^{(V)}_{E_i\cap V}=I_V
\end{equation}
where $I_V$ denotes the identity transformation on $V$.

The equality case of Proposition~\ref{GBLconstanthighdim} indicates why Lemma~\ref{GBLsquarenorm} is important.

\begin{lemma}
\label{GBLsquarenorm}
For linear subspaces $E_1,\ldots,E_k$ of $\R^n$, $n\geq 1$ and $c_1,\ldots,c_k>0$ satisfying
\eqref{highdimcond}, if $\Phi$ is a positive definite linear transform whose eigenspaces are
critical subspaces, then for any $x\in \R^n$, we have
\begin{equation}
\label{GBLsquarenormeq}
\|\Phi x\|^2=
\min_{x=\sum_{i=1}^kc_i  x_i \atop x_i\in E_i}\sum_{i=1}^kc_i\|\Phi x_i\|^2.
\end{equation}
\end{lemma}
\proof  We may assume that ${\rm dim}\,E_i\geq 1$ for $i=1,\ldots,k$.

As the eigenspaces of $\Phi$ are critical subspaces, we deduce that
\begin{equation}
\label{PhiEiEibot}
\Phi(E_i)=E_i\mbox{ \ and \ }\Phi(E_i^\bot)=E_i^\bot.
\end{equation}

For any $x\in\R^n$, we have  $\Phi P_{E_i}x=P_{E_i}\Phi x$ for $i=1,\ldots,k$
by \eqref{PhiEiEibot}; therefore, Lemma~\ref{highdimcritical} (i) yields
\begin{equation}
\label{GBLsquarenormeq0}
\langle \Phi x,\Phi x\rangle=\sum_{i=1}^kc_i\| P_{E_i}\Phi x\|^2= \sum_{i=1}^kc_i\|\Phi P_{E_i}x\|^2.
\end{equation}
Since $x= \sum_{i=1}^kc_i P_{E_i}x$ by \eqref{highdimcond}, we may choose $x_i=P_{E_i}x$
in \eqref{GBLsquarenormeq}, and we have equality in \eqref{GBLsquarenormeq} in this case.
Therefore, Lemma~\ref{GBLsquarenorm} is equivalent to proving that if
$x=\sum_{i=1}^kc_i  x_i$ for $x_i\in E_i$, $i=1,\ldots,k$, then
\begin{equation}
\label{GBLsquarenormineq}
\|\Phi x\|^2\leq \sum_{i=1}^kc_i\|\Phi x_i\|^2.
\end{equation}

\noindent {\bf Case 1 } ${\rm dim}\,E_i=1$ for $i=1,\ldots,k$ and $\Phi=I_n$\\

Let $E_i=\R u_i$ for $u_i\in S^{n-1}$. If $x\in \R^n$, then $P_{E_i}x=\langle u_i,x\rangle u_i$ for $i=1,\ldots,k$,
and \eqref{GBLsquarenormeq0} yields that
$$
\langle x,x\rangle= \sum_{i=1}^kc_i\langle u_i,x\rangle^2.
$$
In addition, any $x_i\in E_i$ is of the form $x_i=t_iu_i$ for $i=1,\ldots,k$ where $\|x_i\|^2=t_i^2$.
If $x=\sum_{i=1}^k c_it_iu_i$, then the H\"older inequality yields
$$
\langle x,x\rangle=\left\langle x,\sum_{i=1}^k c_it_iu_i\right\rangle=
\sum_{i=1}^k c_it_i\langle x,u_i\rangle\leq
\sqrt{\sum_{i=1}^k c_it_i^2}\cdot \sqrt{\sum_{i=1}^k c_i\langle x,u_i\rangle^2}
=\sqrt{\sum_{i=1}^k c_it_i^2}\cdot \sqrt{\langle x,x\rangle},
$$
proving \eqref{GBLsquarenormineq} in this case.\\

\noindent {\bf Case 2 } The general case, $E_1,\ldots,E_k$ and $\Phi$ are as in Lemma~\ref{GBLsquarenorm}\\

Let $V_1,\ldots,V_m$, $m\geq 1$, be the eigenspaces of $\Phi$ corresponding to the eigenvalues
$\lambda_1,\ldots,\lambda_m$.
As $V_1,\ldots,V_m$ are orthogonal critical subspaces and $\R^n=\oplus_{j=1}^mV_j$
As $V_1,\ldots,V_m$ are orthogonal critical subspaces and $\R^n=\oplus_{j=1}^mV_j$, we deduce that
$x_{ij}=P_{V_j}x_i\in E_i\cap V_j$  for any $i=1,\ldots,k$ and $j=1,\ldots,m$, and
$x_i=\sum_{j=1}^m x_{ij}$ for any $i=1,\ldots,k$. It follows that
\begin{eqnarray}
\nonumber
x&=&\sum_{j=1}^m\left( \sum_{E_i\cap V_j\neq\{0\}}c_i x_{ij}\right) \mbox{ \ where}\\
\label{pVjxcixij}
P_{V_j}x&=& \sum_{E_i\cap V_j\neq\{0\}}c_i x_{ij}.
\end{eqnarray}
For any $i=1,\ldots,k$, the vectors $\Phi x_{ij}=\lambda_j x_{ij}$ are pairwise orthogonal for $j=1,\ldots,m$, thus
$$
\sum_{i=1}^kc_i\|\Phi x_i\|^2=\sum_{i=1}^k\left(\sum_{j=1}^m c_i \|\Phi x_{ij}\|^2\right)=
\sum_{j=1}^m\left( \sum_{E_i\cap V_j\neq\{0\}}c_i \|\Phi x_{ij}\|^2\right).
$$
Since $\|\Phi x\|^2=\sum_{j=1}^m \|P_{V_j}\Phi x\|^2=\sum_{j=1}^m \|\Phi P_{V_j}x\|^2$,
\eqref{GBLsquarenormineq} follows if for any $j=1,\ldots,m$, we have
\begin{equation}
\label{GBLsquarenormineqj}
\|\Phi P_{V_j}x\|^2\leq \sum_{E_i\cap V_j\neq\{0\}}c_i\|\Phi x_{ij}\|^2.
\end{equation}

To prove \eqref{GBLsquarenormineqj}, if $E_i\cap V_j\neq\{0\}$, then let ${\rm dim}(E_i\cap V_j)=n_{ij}$,
and let $u_1^{(ij)},\ldots,u_{n_{ij}}^{(ij)}$
be an orthonormal basis  of $E_i\cap V_j$. Since $V_j$ is a critical subspace (see \eqref{highdimcriticalsum}),
if $z\in V_j$, then
 \begin{equation}
\label{VjBLcond}
z=\sum_{i=1}^k c_iP_{E_i} z=\sum_{E_i\cap V_j\neq\{0\}} c_iP_{E_i\cap V_j} z=
\sum_{E_i\cap V_j\neq\{0\}}\sum_{\alpha=1}^{n_{ij}} c_i\langle u_\alpha^{(ij)},z\rangle u_\alpha^{(ij)}.
\end{equation}
\eqref{VjBLcond} shows that the system of all $u_1^{(ij)},\ldots,u_{n_{ij}}^{(ij)}$ when $E_i\cap V_j\neq\{0\}$
form a rank one Brascamp-Lieb data where the coefficient corresponding to $u_\alpha^{(ij)}$ is $c_i$.

According to \eqref{pVjxcixij}, we have
$$
P_{V_j}x= \sum_{E_i\cap V_j\neq\{0\}}\sum_{\alpha=1}^{n_{ij}}c_i \langle u_\alpha^{(ij)},x_{ij}\rangle u_\alpha^{(ij)}.
$$
We deduce from Case~1
applying to $P_{V_j}x$ to the rank one Brascamp-Lieb data  in $V_j$ above that
\begin{eqnarray*}
\|\Phi P_{V_j}x\|^2&=&\lambda_j^2\|P_{V_j}x\|^2\leq \lambda_j^2
\sum_{E_i\cap V_j\neq\{0\}}\sum_{\alpha=1}^{n_{ij}}c_i \langle u_\alpha^{(ij)},x_{ij}\rangle^2\\
&=&
\lambda_j^2\sum_{E_i\cap V_j\neq\{0\}}c_i \|x_{ij}\|^2=
\sum_{E_i\cap V_j\neq\{0\}}c_i\|\Phi x_{ij}\|^2,
\end{eqnarray*}
proving \eqref{GBLsquarenormineqj}, and in turn \eqref{GBLsquarenormineq}
that is equivalent to Lemma~\ref{GBLsquarenorm}.
\proofbox

We now use Proposition~\ref{GBLconstanthighdim} and  Lemma~\ref{GBLsquarenorm} to exhibit the basic type of Gaussian exemizers of Barthe's inequality.

\begin{prop}
\label{GaussianExtremizers}
For linear subspaces $E_1,\ldots,E_k$ of $\R^n$, $n\geq 1$ and $c_1,\ldots,c_k>0$ satisfying
\eqref{highdimcond}, if $\Phi$ is a positive definite linear transform whose eigenspaces are
critical subspaces, then
$$
\int_{\R^n}\left(\sup_{x=\sum_{i=1}^kc_ix_i\atop x_i\in E_i}\;\prod_{i=1}^ke^{-c_i\|\Phi x_i\|^2}\right)dx=
\prod_{i=1}^k\left(\int_{E_i}e^{-\|\Phi x_i\|^2}\,dx_i\right)^{c_i}.
$$
\end{prop}
\proof Let $\widetilde{\Phi}=\pi^{-\frac12}\,\Phi$.
For $i=1,\ldots,k$, let $A_i=\widetilde{\Phi}|_{E_i}$, and hence $A_i:E_i\to E_i$ as the eigenspaces of
$\widetilde{\Phi}$ are critical subspaces.
We deduce first using Lemma~\ref{GBLsquarenorm}, and then the equality case of Proposition~\ref{GBLconstanthighdim} that
\begin{eqnarray*}
\int_{\R^n}\left(\sup_{x=\sum_{i=1}^kc_ix_i\atop x_i\in E_i}\;\prod_{i=1}^ke^{-c_i\|\Phi x_i\|^2}\right)dx
&=&
\int_{\R^n}e^{-\pi\|\widetilde{\Phi} x\|^2}\,dx=\left(\det\widetilde{\Phi}\right)^{-1}=
\prod_{i=1}^k(\det A_i)^{-c_i}\\
&=&\prod_{i=1}^k\left(\int_{E_i}e^{-\pi\|\widetilde{\Phi} x_i\|^2}\,dx_i\right)^{c_i}=
\prod_{i=1}^k\left(\int_{E_i}e^{-\|\Phi x_i\|^2}\,dx_i\right)^{c_i},
\end{eqnarray*}
proving Proposition~\ref{GaussianExtremizers}.
\proofbox

\section{Splitting smooth extremizers along independent and dependent subspaces}
\label{secSplitting}

Optimal transportion as a tool proving geometric inequalities was introduced by Gromov in his Appendix to \cite{MiS86} in the case of the Brunn-Minkowski inequality.
Actually, Barthe's inequality in \cite{Bar98} was one of the first inequalities in probability, analysis or geometry that was obtained via optimal transportation.

We write $\nabla \Theta$ to denote the first derivative of a $C^1$ vector valued function $\Theta$ defined on an open subset of
$\R^n$, and $\nabla^2\varphi$ to denote the Hessian of a real $C^2$ function $\varphi$.
We recall that a vector valued function $\Theta$ on an open set $U\subset\R^n$ is $C^\alpha$ for $\alpha\in (0,1)$
if for any $x_0\in U$ there exist an open neighbourhood $U_0$ of $x_0$ and a $c_0>0$ such that
$\|\Theta(x)-\Theta(y)\|\leq c_0\|x-y\|^\alpha$ for $x,y\in U_0$. In addition, a real function $\varphi$ is $C^{2,\alpha}$
if $\varphi$ is $C^2$ and $\nabla^2 \varphi$ is $C^\alpha$.

Combining Corollary~2.30, Corollary~2.32, Theorem~4.10 and Theorem~4.13 in
Villani \cite{Vil03} on the Brenier map based on McCann \cite{McC95,McC97} for the first two, and on
Caffarelli \cite{Caf90a,Caf90b,Caf92} for the last two theorems, we deduce the following:

\begin{theo}[Brenier, McCann, Caffarelli]
\label{Breniermap}
If $f$ and $g$ are positive $C^\alpha$ probability density functions
 on $\R^n$, $n\geq 1$, for $\alpha\in (0,1)$, then
there exists a  $C^{2,\alpha}$ convex function $\varphi$ on $\R^n$ (unique up to additive constant)
such that $T=\nabla\varphi:\,\R^n\to\R^n$ is bijective and
\begin{equation}
\label{MongeAmpere}
g(x)=f(T(x))\cdot \det \nabla T(x) \mbox{ \ for $x\in\R^n$}.
\end{equation}
\end{theo}
\noindent{\bf Remarks }
The derivative $T=\nabla\varphi$ is the Brenier (transportation) map pushing forward
 the measure on $\R^n$ induced by $g$ to the measure associated to $f$; namely,
$\int_{T(X)}f=\int_X g$ for any measurable $X\subset \R^n$.

In addition, $\nabla T=\nabla^2 \varphi$ is a positive definite symmetrix matrix in Theorem~\ref{Breniermap}, and if $f$ and $g$ are $C^k$ for $k\geq 1$, then $T$ is $C^{k+1}$.

Sometimes it is practical to consider the case $n=0$, when we set $T:\{0\}\to\{0\}$ to be the trivial map.\\

\noindent{\it Proof of Theorem~\ref{RBLtheo} based on Barthe \cite{Bar98}.} First we assume
that each $f_i$ is a $C^1$ positive probability density function
 on $\R^n$, and let us consider the Gaussian densiy $g_i(x)=e^{-\pi\|x\|^2}$ for $x\in E_i$. According to Theorem~\ref{Breniermap},
if $i=1,\ldots,k$, then there exists a
 $C^3$ convex function $\varphi_i$ on $E_i$ such that for the $C^2$ Brenier map $T_i=\nabla\varphi_i$, we have
\begin{equation}
\label{MongeAmperegifi}
g_i(x)=\det \nabla T_i(x) \cdot f_i(T_i(x))\mbox{ \ for all $x\in E_i$}.
\end{equation}
It follows from the Remark after Theorem~\ref{Breniermap} that
$\nabla T_i=\nabla^2 \varphi_i(x)$ is positive definite symmetric matrix for  all $x\in E_i$.
For the $C^2$ transformation $\Theta:\R^n\to\R^n$ given by
\begin{equation}
\label{Thetadef}
\Theta(y)=\sum_{i=1}^kc_iT_i\left(P_{E_i}y\right),\qquad y\in\R^n,
\end{equation}
its differential
$$
\nabla\Theta(y)=\sum_{i=1}^kc_i\nabla T_i\left(P_{E_i}y\right)
$$
is positive definite by Proposition~\ref{GBLconstanthighdim}. It follows that
$\Theta:\R^n\to\R^n$ is
 injective (see \cite{Bar98}), and actually a diffeomorphism. Therefore Proposition~\ref{GBLconstanthighdim}, \eqref{MongeAmperegifi}
and Lemma~\ref{highdimcritical} (i) imply
\begin{align}
\nonumber
&\int_{*,\R^n}\sup_{x=\sum_{i=1}^kc_ix_i,\, x_i\in E_i}\;\prod_{i=1}^kf_i(x_i)^{c_i}\,dx\nonumber\\
&\qquad \geq
\int_{*,\R^n}\left(\sup_{\Theta(y)=\sum_{i=1}^kc_ix_i,\, x_i\in E_i}\;\prod_{i=1}^kf_i(x_i)^{c_i}\right)
\det\left( \nabla\Theta(y)\right)\,dy\nonumber\\
\nonumber
&\qquad\geq  \int_{\R^n}\left(\prod_{i=1}^kf_i\left(T_i\left(P_{E_i}y\right)\right)^{c_i} \right)
\det\left(\sum_{i=1}^kc_i\nabla T_i\left(P_{E_i}y\right)\right)\,dy\\
\label{RBLstep}
&\qquad\geq \int_{\R^n}\left(\prod_{i=1}^kf_i\left( T_i\left(P_{E_i}y\right)\right)^{c_i}\right)
\prod_{i=1}^k\left(\det\nabla T_i\left(P_{E_i}y\right)\right)^{c_i}\,dy\\
\nonumber
&\qquad = \int_{\R^n}\left(\prod_{i=1}^kg_i\left(P_{E_i}y\right)^{c_i}\right)
\,dy=  \int_{\R^n}e^{-\pi \|y\|^2}\,dy=1.
\end{align}
Finally, Barthe's inequality (\ref{RBL}) for arbitrary non-negative integrable functions $f_i$ follows by scaling and approximation (see Barthe \cite{Bar98}).
\proofbox

We now prove that if equality holds in Barthe's inequality (\ref{RBL}), then the diffeomorphism $\Theta$ in \eqref{Thetadef} in the proof of Barthe's inequality splits along the independent subspaces and the dependent subspace.
First we explain how Barthe's inequality behaves under the shifts of the functions involved.
Given proper linear subspaces $E_1,\ldots,E_k$ of $\R^n$ and $c_1,\ldots,c_k>0$ satisfying
\eqref{highdimcond}, first we discuss in what sense Barthe's inequality is translation invariant.
For non-negative integrable function $f_i$ on $E_i$, $i=1,\ldots,k$, let us define
$$
F(x)=\sup_{x=\sum_{i=1}^kc_ix_i,\, x_i\in E_i}\;\prod_{i=1}^kf_i(x_i)^{c_i}.
$$
We observe that for any $e_i\in E_i$, defining $\tilde{f}_i(x)=f_i(x+e_i)$ for $x\in E_i$, $i=1,\ldots,k$, we have
\begin{equation}
\label{RBLtranslation}
\widetilde{F}(x)=\sup_{x=\sum_{i=1}^kc_ix_i,\, x_i\in E_i}\;\prod_{i=1}^k\tilde{f}_i(x_i)^{c_i}=
F\left(x+\sum_{i=1}^kc_ie_i\right).
\end{equation}

\begin{prop}
\label{Thetasplits}
For  non-trivial linear subspaces  $E_1,\ldots,E_k$ of $\R^n$ and $c_1,\ldots,c_k>0$ satisfying
\eqref{highdimcond0}, we write $F_1,\ldots,F_l$ to denote the independent subspaces (if exist), and
$F_0$ to denote the dependent subspace (possibly $F_0=\{0\}$).
Let us assume that equality holds in \eqref{RBL}
for positive $C^1$ probability densities $f_i$ on $E_i$, $i=1,\ldots,k$, let
$g_i(x)=e^{-\pi\|x\|^2}$ for $x\in E_i$, let $T_i:E_i\to E_i$ be the $C^2$ Brenier map satisfying
\begin{equation}
\label{MongeAmperegifi0}
g_i(x)=\det \nabla T_i(x) \cdot f_i(T_i(x))\mbox{ \ for all $x\in E_i$},
\end{equation}
and let
$$
\Theta(y)=\sum_{i=1}^kc_iT_i\left(P_{E_i}y\right),\qquad y\in\R^n.
$$
\begin{description}
\item{(i)} For any $i\in\{1,\ldots,k\}$ there exists positive $C^1$
integrable $h_{i0}:\,F_0\cap E_i\to[0,\infty)$ (where $h_{i0}(0)=1$ if $F_0\cap E_i=\{0\}$),
and for any $i\in\{1,\ldots,k\}$ and $j\in\{1,\ldots,l\}$ with $F_j\subset E_i$,
 there exists positive $C^1$
integrable $h_{ij}:\,F_j\to[0,\infty)$  such that
$$
f_i(x)=h_{i0}(P_{F_0}x)\cdot
\prod_{F_j\subset E_i\atop j\geq 1}h_{ij}(P_{F_j}x) \mbox{ \ \ \  for $x\in E_i$}.
$$
\item{(ii)} For $i=1,\ldots,k$, $T_i(E_i\cap F_p)=E_i\cap F_p$ whenever $E_i\cap F_p\neq \{0\}$ for $p\{0,\ldots,l\}$, and
 if $x\in E_i$, then
$$
T_i(x)=\bigoplus_{E_i\cap F_p\neq \{0\} \atop p\geq 0} T_i(P_{F_p}x).
$$

\item{(iii)} For $i=1,\ldots,k$, there exist $C^2$ functions $\Omega_i:E_i\to E_i$
and $\Gamma_i:E_i^\bot\to E_i^\bot$ such that
$$
\Theta(y)=\Omega_i(P_{E_i}y)+\Gamma_i(P_{E_i^\bot}y) \mbox{ \ \ \  for $y\in \R^n$}.
$$
\item{(iv)} If $y\in\R^n$, then the eigenspaces of the positive definite matrix $\nabla\Theta(y)$ are critical subspaces,
and $\nabla T_i(P_{E_i}y)=\nabla\Theta(y)|_{E_i}$ for $i=1,\ldots,k$.
\end{description}
\end{prop}
\proof
According to \eqref{RBLtranslation}, we may assume that
\begin{equation}
\label{Tizero}
T_i(0)=0\mbox{ \ \ for $i=1,\ldots,k$},
\end{equation}

If equality holds in \eqref{RBL}, then equality holds in the determinantal inequality in \eqref{RBLstep} in the proof of
Barthe's inequality; therefore,
we apply the equality case of Proposition~\ref{GBLconstanthighdim}.
In particular, for any $x\in\R^n$,
there exist $m_x\geq 1$ and linear subspaces $V_{1,x},\ldots,V_{m_x,x}$ where $V_{1,x}=\R^n$ if $m_x=1$, and
$V_{1,x},\ldots,V_{m_x,x}$ are pairwise orthogonal indecomposable critical subspaces spanning $\R^n$
if $m_x\geq 2$, and
there exist $\lambda_{1,x},\ldots,\lambda_{m_x,x}>0$ such that
if $E_i\cap V_{j,x}\neq\{0\}$, then
\begin{equation}
\label{nablaTcond}
\nabla T_i(P_{E_i}x)|_{E_i\cap V_{j,x}}=\lambda_{j,x}I_{E_i\cap V_{j,x}};
\end{equation}
and in addition, each $E_i$ satisfies ({\it cf.} \eqref{highdimcriticalV1V2})
\begin{equation}
\label{VjxEi}
E_i=\oplus_{E_i\cap V_{j,x}\neq\{0\}}E_i\cap V_{j,x}.
\end{equation}

Let us consider a fixed $E_i$, $i\in\{1,\ldots,k\}$.
First we claim that if $y\in E_i$, then
\begin{equation}
\label{nablaTiEicapFp}
\begin{array}{rcll}
\nabla T_i(y)(F_p)&=&F_p&\mbox{ \ if $p\geq 1$ and $E_i\cap F_p\neq\{0\}$}\\[1ex]
\nabla T_i(y)(F_0\cap E_i)&=&F_0\cap E_i.&
\end{array}
\end{equation}
To prove \eqref{nablaTiEicapFp}, we take $y=x$ in \eqref{nablaTcond}.
 If $p\geq 1$ and $E_i\cap F_p\neq\{0\}$, then $F_p\subset E_i$, and
 Lemma~\ref{indec-dependent-critical} yields that
\begin{eqnarray*}
\oplus_{F_p\cap V_{j,y}\neq \{0\}}  V_{j,y}&\subset & F_p\\
\oplus_{F_p\cap V_{j,y}= \{0\}}  V_{j,y}&\subset & F_p^\bot.
\end{eqnarray*}
Since the subspaces $V_{j,y}$ span $\R^n$, we have
$$
F_p=\oplus_{E_i\cap V_{j,y}\neq \{0\} \atop V_{j,y}\subset F_p}  V_{j,y};
$$
therefore, \eqref{nablaTcond} implies \eqref{nablaTiEicapFp} if $p\geq 1$.

For the case of $F_0$ in \eqref{nablaTiEicapFp}, it follows from \eqref{VjxEi} and Lemma~\ref{indec-dependent-critical}  that
if $E_i\cap F_0\neq \{0\}$, then
\begin{equation}
\label{VjxEiF0}
E_i\cap F_0=\oplus_{E_i\cap V_{j,y}\neq\{0\}\atop
 V_{j,y}\subset F_0}E_i\cap V_{j,y}.
\end{equation}
Therefore, \eqref{nablaTcond} completes the proof of \eqref{nablaTiEicapFp}.

It follows from \eqref{nablaTiEicapFp} that if $E_i\cap F_p\neq\{0\}$, $y\in E_i$, $v\in E_i\cap F_p\cap S^{n-1}$ and
$w\in E_i\cap F_p^\bot\cap S^{n-1}$, then
\begin{equation}
\label{nablaTisumFp}
\left\langle v, \left.\frac{\partial}{\partial t}T_i(y+tw)\right|_{t=0}\right\rangle=0.
\end{equation}
In turn, \eqref{nablaTiEicapFp}, \eqref{nablaTisumFp} and $T_i(0)=0$ ({\it cf.} \eqref{Tizero}) imply that
if $y\in E_i$, then
\begin{eqnarray}
\label{TisumFp}
T_i(E_i\cap F_p)&=&E_i\cap F_p \mbox{ \ whenever $E_i\cap F_p\neq \{0\}$ for $p\geq 0$,}\\
\label{TisumFpy}
T_i(y)&=&\bigoplus_{E_i\cap F_p\neq \{0\} \atop p\geq 0} T_i(P_{F_p}y).
\end{eqnarray}
We deduce from \eqref{TisumFpy} that if $y\in E_i$, then
\begin{eqnarray}
\label{nablaTisumFpdet}
\det\nabla T_i(y)&=&\prod_{E_i\cap F_p\neq \{0\}\atop p\geq 0 } \det\left(\nabla T_i(P_{F_p}y)|_{F_p}\right).
\end{eqnarray}
We conclude (i) from \eqref{nablaTisumFp}, \eqref{TisumFp}, \eqref{TisumFpy},  
and  \eqref{nablaTisumFpdet}
as \eqref{MongeAmperegifi0} yields that if $y\in E_i$, then
$$
 f_i(T_i(y))=\prod_{E_i\cap F_p\neq \{0\}\atop p\geq 0 } \frac{e^{-\pi\|P_{F_p}y\|^2}}
{\det\left(\nabla T_i(P_{F_p}y)|_{F_p}\right)}.
$$

We deduce (ii) from \eqref{TisumFp} and \eqref{TisumFpy}.

For (iii), it follows from Proposition~\ref{GBLconstanthighdim} that for any $x\in\R^n$, the spaces $V_{j,x}$ are
eigenspaces for $\nabla\Theta(x)$ and span $\R^n$; therefore, \eqref{highdimcriticalV} implies that if
$x\in\R^n$ and $i\in\{1,\ldots,k\}$, then
$$
\nabla \Theta(x)=\nabla \Theta(x)|_{E_i}\oplus \nabla \Theta(x)|_{E_i^\bot}.
$$
Since $\Theta(0)=0$ by \eqref{Tizero},  for fixed $i\in\{1,\ldots,k\}$, we conclude
\begin{eqnarray*}
\Theta(E_i)&=&E_i;\\
\Theta(x)&=&\left.\Theta\left(P_{E_i}x\right)\right|_{E_i}\oplus \left.\Theta\left(P_{E_i^\bot}x\right)\right|_{E_i^\bot}
\mbox{ \ if $x\in\R^n$}.
\end{eqnarray*}

Finally, (iv) directly follows from Proposition~\ref{GBLconstanthighdim}, completing the proof of Proposition~\ref{Thetasplits}.
\proofbox

Next we show that if the extremizers $f_1,\ldots,f_k$ in Proposition~\ref{Thetasplits} are of the form as in (i), then
for any given $F_j\neq\{0\}$, the functions
$h_{ij}$ on $F_j$ for all $i$ with $E_i\cap F_j\neq\{0\}$ are also extremizers.
We also need the Pr\'ekopa-Leindler inequality
Theorem~\ref{PL}
(proved in various forms by Pr\'ekopa \cite{Pre71,Pre73}, Leindler \cite{Lei72} and
 Borell \cite{Bor75}) whose equality case was clarified by Dubuc \cite{Dub77}
(see the survey Gardner \cite{gardner}). In turn, the Pr\'ekopa-Leindler inequality
\eqref{PLineq} is of the very similar structure
like Barthe's inequality \eqref{RBL}. Again, the inequality is usually stated using outer integrals, but being the special case of Barthe's inequality \eqref{RBL}, the remarks after Theorem~\ref{RBLtheo} concerning inner integration apply.

\begin{theo}[Pr\'ekopa, Leindler, Dubuc]
\label{PL}
For $m\geq 2$, $\lambda_1,\ldots,\lambda_m\in(0,1)$
with $\lambda_1+\ldots+\lambda_m=1$ and
 integrable $\varphi_1,\ldots,\varphi_m:\,\R^n\to[0,\infty)$, we have
\begin{equation}
\label{PLineq}
\int_{*,\R^n} \sup_{x=\sum_{i=1}^m\lambda_ix_i,\, x_i\in \R^n}\;\prod_{i=1}^m\varphi_i(x_i)^{\lambda_i}\,dx
\geq \prod_{i=1}^m\left(\int_{\R^n}\varphi_i\right)^{\lambda_i},
\end{equation}
and if equality holds and the left hand side is positive and finite, then there exist
a log-concave function $\varphi$ and
 $a_i>0$ and $b_i\in\R^n$ for $i=1,\ldots,m$ such that
$$
\varphi_i(x)=a_i\, \varphi(x-b_i)
$$
for Lebesgue a.e. $x\in\R^n$, $i=1,\ldots,m$.
\end{theo}

For linear subspaces  $E_1,\ldots,E_k$ of $\R^n$ and $c_1,\ldots,c_k>0$ satisfying
\eqref{highdimcond0}, we assume that $F_{\rm dep}\neq \R^n$, and write $F_1,\ldots,F_l$ to denote the independent subspaces. We verify that if $j\in\{1,\ldots,l\}$, then
\begin{equation}
\label{sumci1}
\sum_{E_i\supset F_j}c_i=1.
\end{equation}
For this, let $x\in F_j\backslash\{0\}$. We observe that for any $E_i$,
either $F_j\subset E_i$, and hence $P_{E_i}x=x$, or $F_j\subset E_i^\bot$,
and hence $P_{E_i}x=o$. We deduce from \eqref{highdimcond0} that
$$
x=\sum_{i=1}^k c_iP_{E_i}x=\left(\sum_{F_j\subset E_i} c_i\right)\cdot x,
$$
which formula in turn implies \eqref{sumci1}.

\begin{prop}
\label{Thetasplitshij}
For linear subspaces  $E_1,\ldots,E_k$ of $\R^n$ and $c_1,\ldots,c_k>0$ satisfying
\eqref{highdimcond0}, we write $F_1,\ldots,F_l$ to denote the independent subspaces (if exist), and
$F_0$ denote the dependent subspace (possibly $F_0=\{0\}$).
Let us assume that equality holds in Barthe's inequality \eqref{RBL}
for probability densities $f_i$ on $E_i$, $i=1,\ldots,k$, and
for any $i\in\{1,\ldots,k\}$ there exists
integrable $h_{i0}:\,F_0\cap E_i\to[0,\infty)$ (where $h_{i0}(0)=1$ if $F_0\cap E_i=\{0\}$),
and for any $i\in\{1,\ldots,k\}$ and $j\in\{1,\ldots,l\}$ with $F_j\subset E_i$,
 there exists non-negative
integrable $h_{ij}:\,F_j\to[0,\infty)$  such that
\begin{equation}
\label{fihijcond}
f_i(x)=h_{i0}(P_{F_0}x)\cdot
\prod_{F_j\subset E_i\atop j\geq 1}h_{ij}(P_{F_j}x) \mbox{ \ \ \  for $x\in E_i$}.
\end{equation}
\begin{description}
\item {(i)} If $F_0\neq \{0\}$, then
$\sum_{E_i\cap F_0\neq\{0\}}c_iP_{E_i\cap F_0}={\rm Id}_{F_0}$ and
$$
\int_{*,F_0}\sup_{x=\sum\{c_ix_i:\, x_i\in E_i\cap F_0\,\&\,E_i\cap F_0\neq\{0\}\}}\;
\prod_{E_i\cap F_0\neq\{0\}}h_{i0}(x_i)^{c_i}\,dx
= \prod_{E_i\cap F_0\neq\{0\}}\left(\int_{E_i\cap F_0}h_{i0}\right)^{c_i}.
$$
\item {(ii)} If $F_0\neq \R^n$, then there exist integrable $\psi_{j}:\,F_j\to[0,\infty)$
for $j=1,\ldots,l$ where $\psi_j$ is log-concave whenever $F_j\subset E_\alpha\cap E_\beta$ for $\alpha\neq \beta$, and there exist $a_{ij}>0$ and $b_{ij}\in F_j$
for any $i\in\{1,\ldots,k\}$ and $j\in\{1,\ldots,l\}$ with $F_j\subset E_i$ such that
$h_{ij}(x)=a_{ij}\cdot \psi_j(x-b_{ij})$ for $i\in\{1,\ldots,k\}$ and $j\in\{1,\ldots,l\}$ with $F_j\subset E_i$.
\end{description}
\end{prop}
\proof We only present the argument in the case $F_0\neq \R^n$ and $F_0\neq \{0\}$. If
$F_0= \R^n$, then the same argument works ignoring the parts involving $F_1,\ldots,F_l$, and if
$F_0= \{0\}$, then the same argument works ignoring the parts involving $F_0$.

Since $F_0\oplus F_1\oplus\ldots \oplus F_l=\R^n$ and
$F_0,\ldots, F_l$ are critical subspaces,
\eqref{highdimcriticalV1V2} yields for $i=1,\ldots,k$ that
\begin{equation}
\label{EiFjsum}
E_i=(E_i\cap F_0)\oplus \bigoplus_{F_j\subset E_i\atop j\geq 1} F_j;\end{equation}
therefore, the Fubini theorem and \eqref{fihijcond} imply that
\begin{equation}
\label{fihijprod1}
\int_{E_i}f_i=\left(\int_{E_i\cap  F_0}h_{i0}\right)\cdot
\prod_{F_j\subset E_i\atop j\geq 1}\int_{F_j}h_{ij}.
\end{equation}
On the other hand, using again $ F_0\oplus F_1\oplus\ldots \oplus F_l=\R^n$, we deduce that
if $x=\sum_{j=0}^l z_j$ where  $z_j\in F_j$ for $j\geq 0$,
then $z_j=P_{F_j}x$.
It follows from \eqref{EiFjsum} that for any $x\in\R^n$, we have
\begin{eqnarray*}
\sup_{x=\sum_{i=1}^kc_ix_i,\atop x_i\in E_i}\;\prod_{i=1}^kf_i(x_i)^{c_i}&=&
\left( \sup_{P_{ F_0}x=\sum_{i=1}^kc_ix_{0i},\atop x_{0i}\in E_i\cap F_0}\;
\prod_{i=1}^kh_{i0}(x_{i0})\right)\times\\
&&\times\prod_{j=1}^l
\left( \sup_{P_{ F_j}x=\sum_{F_j\subset E_i}c_ix_{ji},\atop x_{ji}\in F_j}
\prod_{F_j\subset E_i}h_{ij}(x_{ji})^{c_i}\right),
\end{eqnarray*}
and hence
\begin{eqnarray}
\label{fihijprod2}
\int_{*,\R^n}\sup_{x=\sum_{i=1}^kc_ix_i,\atop x_i\in E_i}\;\prod_{i=1}^kf_i(x_i)^{c_i}\,dx&=&
\left( \int_{*,F_0}\sup_{x=\sum_{i=1}^kc_ix_{i},\atop x_{i}\in E_i\cap F_0}\;
\prod_{i=1}^kh_{i0}(x_{i})\,dx\right)\times\\
\nonumber
&&\times \prod_{j=1}^l
\left( \int_{*,F_j}\sup_{x=\sum_{F_j\subset E_i}c_ix_{i},\atop x_{i}\in F_j}
\prod_{F_j\subset E_i}h_{ij}(x_{i})^{c_i}\,dx\right).
\end{eqnarray}

As $F_0$ is a critical subspace, we have
$$
\sum_{i=1}^kc_iP_{E_i\cap F_0}={\rm Id}_{F_0},
$$
and hence Barthe's inequality \eqref{RBL} yields
\begin{equation}
\label{fihijprod3}
\int_{*,F_0}\sup_{x=\sum_{i=1}^kc_ix_{i},\atop x_{i}\in E_i\cap F_0}\;
\prod_{i=1}^kh_{i0}(x_{i})\,dx\geq
\prod_{i=1}^k\left(\int_{E_i\cap  F_0}h_{i0}\right)^{c_i}.
\end{equation}

We deduce from \eqref{sumci1} and the Pr\'ekopa-Leindler inequality \eqref{PLineq} that
if $j=1,\ldots,l$, then
\begin{equation}
\label{fihijprod4}
 \int_{*,F_j}\sup_{x=\sum_{F_j\subset E_i}c_ix_{i},\atop x_{i}\in F_j}
\prod_{F_j\subset E_i}h_{ij}(x_{i})^{c_i}\,dx\geq\prod_{E_i\supset F_j}\left(\int_{F_j}h_{ij}\right)^{c_i}.
\end{equation}

Combining \eqref{fihijprod1}, \eqref{fihijprod2}, \eqref{fihijprod3} and \eqref{fihijprod4} with the fact that
$f_1,\ldots,f_k$ are extremizers for Barthe's inequality \eqref{RBL} implies that
if $j=1,\ldots,l$, then
\begin{eqnarray}
\label{fihijprod5}
\int_{*,F_0}\sup_{x=\sum_{i=1}^kc_ix_{i},\atop x_{i}\in E_i\cap F_0}\;
\prod_{i=1}^kh_{i0}(x_{i})\,dx&=&
\prod_{i=1}^k\left(\int_{E_i\cap  F_0}h_{i0}\right)^{c_i}\\
\label{fihijprod6}
 \int_{*,F_j}\sup_{x=\sum_{F_j\subset E_i}c_ix_{i},\atop x_{i}\in F_j}
\prod_{F_j\subset E_i}h_{ij}(x_{i})^{c_i}\,dx&=&\prod_{E_i\supset F_j}\left(\int_{F_j}h_{ij}\right)^{c_i}.
\end{eqnarray}
We observe that \eqref{fihijprod5} is just (i). In addition, (ii) follows from the equality conditions in the Pr\'ekopa-Leindler inequality (see Theorem~\ref{PL}).
\proofbox

\section{Convolution and product of extremizers}
\label{secExtremizers}

Given proper linear subspaces $E_1,\ldots,E_k$ of $\R^n$ and $c_1,\ldots,c_k>0$ satisfying
\eqref{highdimcond}, we say that the non-negative integrable functions $f_1,\ldots,f_k$ with positive integrals are extremizers
if equality holds in \eqref{RBL}. In order to deal with positive smooth functions, we use convolutions. More precisely, Lemma~2 in Barthe \cite{Bar98} states the following.

\begin{lemma}
\label{convolution}
Given proper linear subspaces $E_1,\ldots,E_k$ of $\R^n$ and $c_1,\ldots,c_k>0$ satisfying
\eqref{highdimcond}, if $f_1,\ldots,f_k$ and $g_1,\ldots,g_k$ are extremizers in  Barthe's  inequality \eqref{RBL}, then $f_1*g_1,\ldots,f_k*g_k$
are also are extremizers.
\end{lemma}
\proof We may assume that $\int_{\R^n} f_i=\int_{\R^n} g_i=1$ for $i=1,\ldots,k$. We  define
\begin{eqnarray*}
F(x)&=&\sup_{x=\sum_{i=1}^kc_ix_i,\, x_i\in E_i}\;\prod_{i=1}^kf_i(x_i)^{c_i}\\
G(y)&=&\sup_{y=\sum_{i=1}^kc_iy_i,\, y_i\in E_i}\;\prod_{i=1}^kg_i(y_i)^{c_i}.
\end{eqnarray*}
 Possibly $F$ and $G$ are not measurable but
as $f_1,\ldots,f_k$ and $g_1,\ldots,g_k$ are extremizers, 
we have $\int_{*,\R^n} F=\int_{*,\R^n}G=1$, and 
there exist measurable $\widetilde{F}\geq F$
and $0\leq \widetilde{G}\leq G$ such that
$\int_{\R^n}\widetilde{F}(x)\,dx=\int_{\R^n}\widetilde{G}(x)\,dx=1$, and hence neither $\{F>\widetilde{F}\}$
nor $\{G>\widetilde{G}\}$ contains a subset of $\R^n$ with positive measure. 
We write any point of $\R^{2n}$ in the form $(x,y)$ for $x,y\in\R^n$, and hence
$(z,y)\mapsto \widetilde{F}(z)\widetilde{G}(y)$ is a measurable witness for $(z,y)\mapsto F(z)G(y)$, and in turn
$(x,y)\mapsto \widetilde{F}(x-y)\widetilde{G}(y)$ is a measurable witness for $(x,y)\mapsto F(x-y)G(y)$ in the case of inner integrals.
We deduce writing $x_i=z_i+y_i$ in (\ref{xiyizi}) for $i=1,\ldots,k$ and using Barthe's inequality in (\ref{figiconvRBL}) that
\begin{eqnarray}
\nonumber
1&=& \int_{\R^n} \widetilde{F}*\widetilde{G}(x) dx=\int_{\R^n}  \int_{\R^n} \widetilde{F}(x-y)\widetilde{G}(y) \,dy dx\\
\nonumber
 &=& \int_{*,\R^{2n}}
\sup_{x-y=\sum_{i=1}^kc_iz_i,\, z_i\in E_i}\;\prod_{i=1}^kf_i(z_i)^{c_i}
\sup_{y=\sum_{i=1}^kc_iy_i,\, y_i\in E_i}\;\prod_{i=1}^kg_i(y_i)^{c_i} \,d(x,y)\\
\label{xiyizi}
&=& \int_{*,\R^{2n}}
\sup_{x-y=\sum_{i=1}^kc_iz_i,\, z_i\in E_i}\;
\sup_{y=\sum_{i=1}^kc_iy_i,\, y_i\in E_i}\;
\prod_{i=1}^kf_i(z_i)^{c_i}
\prod_{i=1}^kg_i(y_i)^{c_i} \,d(x,y)\\
\nonumber
&=&\int_{*,\R^{2n}}
\sup_{x=\sum_{i=1}^kc_ix_i,\, x_i\in E_i}\;
\sup_{y=\sum_{i=1}^kc_iy_i,\, y_i\in E_i}\;
\prod_{i=1}^kf_i(x_i-y_i)^{c_i}
\prod_{i=1}^kg_i(y_i)^{c_i} \,d(x,y)\\
\label{figiconvRBL}
&\geq& \int_{*,\R^n}  \sup_{x=\sum_{i=1}^kc_ix_i,\, x_i\in E_i}\;\int_{*,\R^n}
\sup_{y=\sum_{i=1}^kc_iy_i,\, y_i\in E_i}\;
\prod_{i=1}^k\big(f_i(x_i-y_i)g_i(y_i)\big)^{c_i} \,dy dx\\
\nonumber
&\geq& \int_{*,\R^n}  \sup_{x=\sum_{i=1}^kc_ix_i,\, x_i\in E_i}\;\prod_{i=1}^k
\left(\int_{E_i}
f_i(x_i-y_i) g_i(y_i) \,dy_i\right)^{c_i} dx\\
\nonumber
&=& \int_{*,\R^n}  \sup_{x=\sum_{i=1}^kc_ix_i,\, x_i\in E_i}\;\prod_{i=1}^k
\big(f_i*g_i(x_i)\big)^{c_i} dx\geq \prod_{i=1}^k\left(\int_{\R^n}f_i*g_i(x_i)\right)^{c_i}=1
\end{eqnarray}
as $\int_{E_i} f_i*g_i=1$ for $i=1,\ldots,k$. In turn,
we conclude that $f_i*g_i$, $i=1,\ldots,k$, is also an extremizer.
\proofbox

 Since in a certain case we want to work with Lebesgue integral instead of outer integrals, we use the following statement that can be proved via compactness argument.

\begin{lemma}
\label{noouterint}
Given proper linear subspaces $E_1,\ldots,E_k$ of $\R^n$ and $c_1,\ldots,c_k>0$ satisfying
\eqref{highdimcond}, if $h_i$ is a positive continuous functions satisfying
$\lim_{x\to \infty}h_i(x)=0$ for $i=1,\ldots,k$,
 then the function
$$
h(x)=\sup_{x=\sum_{i=1}^kc_ix_i,\atop x_i\in E_i}\prod_{i=1}^kh_i(x_i)^{c_i}
$$
of $x\in\R^n$ is continuous.
\end{lemma}

Next we show that the product of a shift of a smooth extremizer and a Gaussian is also an extremizer
for  Barthe's  inequality.

\begin{lemma}
\label{RBLextproduct}
Given proper linear subspaces $E_1,\ldots,E_k$ of $\R^n$ and $c_1,\ldots,c_k>0$ satisfying
\eqref{highdimcond}, if $f_1,\ldots,f_k$  are positive bounded $C^1$ are extremizers in Barthe's  inequality \eqref{RBL}, and $g_i(x)= e^{-\pi\|x\|^2}$ for $x\in E_i$,
 then there
exist $z_i\in E_i$, $i=1,\ldots,k$, such that the functions $y\mapsto f_i(y-z_i)g_i(y)$ of $y\in E_i$, $i=1,\ldots,k$, are also
extremizers for \eqref{RBL}.
\end{lemma}
\proof We may assume that $f_1,\ldots,f_k$  are probability densities.

Readily the functions $\tilde{f}_1,\ldots,\tilde{f}_k$ defined by
 $\tilde{f}_i(y)=f_i(-y)$ for $y\in E_i$ and
$i=1,\ldots,k$ are also extremizers. We deduce from Lemma~\ref{convolution} that
the functions
$\tilde{f}_i* g_i$ for $i=1,\ldots,k$ are also extremizers
where each $\tilde{f}_i* g_i$ is a probability density on $E_i$.
According to Theorem~\ref{Breniermap},
if $i=1,\ldots,k$, then there exists a
 $C^2$ Brenier map $S_i:E_i\to E_i$ such that
$$
g_i(x)=\det \nabla S_i(x) \cdot (\tilde{f}_i* g_i)(S_i(x))\mbox{ \ for all $x\in E_i$},
$$
and $\nabla S_i(x)$ is a positive definite symmetric matrix for  all $x\in E_i$.
As in the proof of Theorem~\ref{RBLtheo} above, we consider the $C^2$ diffeomorphism $\Theta:\R^n\to\R^n$ given by
$$
\Theta(y)=\sum_{i=1}^kc_iS_i\left(P_{E_i}y\right),\qquad y\in\R^n.
$$
whose positive definite differential is
$$
\nabla\Theta(y)=\sum_{i=1}^kc_i\nabla S_i\left(P_{E_i}y\right).
$$

On the one hand, we note that
if $x=\sum_{i=1}^kc_ix_i$ for $x_i\in E_i$, then
$$
\|x\|^2\leq \sum_{i=1}^kc_i\|x_i\|^2
$$
holds according to Barthe \cite{Bar98}; or equivalently,
$$
\prod_{i=1}^kg_i(x_i)^{c_i}\leq  e^{-\pi\|x\|^2}.
$$
Since $f_i$ is positive, bounded, continuous and in $L_1(E_i)$ for $i=1,\ldots,k$, we observe that the function
\begin{equation}
\label{figiwtendsinfinity}
z\mapsto \int_{\R^n}\sup_{x=\sum_{i=1}^kc_ix_i,\atop x_i\in E_i}\;
\prod_{i=1}^kf_i\left(x_i-S_i(P_{E_i}\Theta^{-1}z)\right)^{c_i}g_i(x_i)^{c_i}\,dx
\end{equation}
of $z\in \R^n$ is continuous.

Using also that $\tilde{f}_1,\ldots,\tilde{f}_k$ are extremizers and probability density functions, we have
\begin{align*}
\int_{*,\R^n}\int_{*,\R^n}\sup_{z=\sum_{i=1}^kc_iz_i,\atop z_i\in E_i}\sup_{x=\sum_{i=1}^kc_ix_i,\atop x_i\in E_i}\;\prod_{i=1}^kf_i(x_i-z_i)^{c_i}g_i(x_i)^{c_i}\,dx\,dz&\\
=\int_{*,\R^n}\int_{*,\R^n}
\sup_{x=\sum_{i=1}^kc_ix_i,\atop x_i\in E_i}
\left(\prod_{i=1}^kg_i(x_i)^{c_i}\right)
\sup_{z=\sum_{i=1}^kc_iz_i,\atop z_i\in E_i}\;\prod_{i=1}^kf_i(x_i-z_i)^{c_i}\,dz\,dx&\\
\leq \int_{*,\R^n} e^{-\pi\|x\|^2}\int_{*,\R^n}
\sup_{x=\sum_{i=1}^kc_ix_i,\atop x_i\in E_i}
\sup_{z=\sum_{i=1}^kc_iz_i,\atop z_i\in E_i}\;\prod_{i=1}^kf_i(x_i-z_i)^{c_i}\,dz\,dx&\\
= \int_{*,\R^n}e^{-\pi\|x\|^2}\int_{*,\R^n}
\sup_{z-x=\sum_{i=1}^kc_iy_i,\atop y_i\in E_i}
\;\prod_{i=1}^k\tilde{f}_i(y_i)^{c_i}\,dz\,dx&\\
= \int_{*,\R^n}e^{-\pi\|x\|^2}\int_{*,\R^n}
\sup_{w=\sum_{i=1}^kc_iy_i,\atop y_i\in E_i}
\;\prod_{i=1}^k\tilde{f}_i(y_i)^{c_i}\,dw\,dx&\\
= \int_{\R^n}e^{-\pi\|x\|^2}\,dx=1.&
\end{align*}

Using
Lemma~\ref{noouterint} and \eqref{figiwtendsinfinity} in \eqref{stillouter},
 Barthe's  inequality \eqref{RBL} in \eqref{RBLfigi}
 and Proposition~\ref{GBLconstanthighdim} in \eqref{ThetaSi}, we deduce that
\begin{align}
\nonumber
1&\geq \int_{*,\R^n}\int_{*,\R^n}\sup_{z=\sum_{i=1}^kc_iz_i\atop z_i\in E_i}\sup_{x=\sum_{i=1}^kc_ix_i\atop x_i\in E_i}\;\prod_{i=1}^kf_i(x_i-z_i)^{c_i}g_i(x_i)^{c_i}\,dx\,dz\\
\label{stillouter}
&\geq
\int_{*,\R^n}\int_{*,\R^n}\sup_{x=\sum_{i=1}^kc_ix_i\atop x_i\in E_i}\;
\prod_{i=1}^kf_i\left(x_i-S_i(P_{E_i}\Theta^{-1}z)\right)^{c_i}g_i(x_i)^{c_i}\,dx\,dz\\
\label{RBLfigi}
&=\int_{\R^n}\int_{\R^n}\sup_{x=\sum_{i=1}^kc_ix_i,\atop x_i\in E_i}\;
\prod_{i=1}^kf_i\left(x_i-S_i(P_{E_i}\Theta^{-1}z)\right)^{c_i}g_i(x_i)^{c_i}\,dx\,dz\\
\nonumber
&\geq \int_{\R^n}\prod_{i=1}^k\left(\int_{E_i}
f_i\left(x_i-S_i(P_{E_i}\Theta^{-1}z)\right)g_i(x_i)\,dx_i\right)^{c_i}\,dz\\
\nonumber
&=\int_{\R^n}\prod_{i=1}^k(\tilde{f}_i* g_i)\left(S_i(P_{E_i}\Theta^{-1}z)\right)^{c_i}\,dz\\
\nonumber
&=
\int_{\R^n} \left(\prod_{i=1}^k(\tilde{f}_i* g_i)\left(S_i\left(P_{E_i}y\right)\right)^{c_i} \right)
\det\left( \nabla\Theta(y)\right)\,dy\\
\label{ThetaSi}
&=  \int_{\R^n}\left(\prod_{i=1}^k(\tilde{f}_i* g_i)\left(S_i\left(P_{E_i}y\right)\right)^{c_i} \right)
\det\left(\sum_{i=1}^kc_i\nabla S_i\left(P_{E_i}y\right)\right)\,dy\\
\nonumber
&
\geq \int_{\R^n}\left(\prod_{i=1}^k(\tilde{f}_i* g_i)\left( S_i\left(P_{E_i}y\right)\right)^{c_i}\right)
\prod_{i=1}^k\left(\det\nabla S_i\left(P_{E_i}y\right)\right)^{c_i}\,dy\\
\nonumber
&
 = \int_{\R^n}\left(\prod_{i=1}^kg_i\left(P_{E_i}y\right)^{c_i}\right)
\,dy=  \int_{\R^n} e^{-\pi\|y\|^2}\,dy=1.
\end{align}
In particular, we conclude that
\begin{eqnarray*}
1&\geq&
\int_{\R^n}\int_{\R^n}\sup_{x=\sum_{i=1}^kc_ix_i,\atop x_i\in E_i}\;\prod_{i=1}^k
f_i\left(x_i-S_i(P_{E_i}\Theta^{-1}z)\right)^{c_i}g_i(x_i)^{c_i}\,dx\,dz\\
&\geq & \int_{\R^n}
\prod_{i=1}^k\left(\int_{E_i}f_i\left(x_i-S_i(P_{E_i}\Theta^{-1}z)\right)g_i(x_i)\,dx_i\right)^{c_i}\,dz\geq 1.
\end{eqnarray*}
Because of  Barthe's  inequality \eqref{RBL}, it follows from \eqref{figiwtendsinfinity}
that
\begin{eqnarray*}
&&\int_{\R^n}\sup_{x=\sum_{i=1}^kc_ix_i,\atop x_i\in E_i}\;\prod_{i=1}^k
f_i\left(x_i-S_i(P_{E_i}\Theta^{-1}z)\right)^{c_i}g_i(x_i)^{c_i}\,dx\\
&=&
\prod_{i=1}^k\left(\int_{E_i}f_i\left(x_i-S_i(P_{E_i}\Theta^{-1}z)\right)g_i(x_i)\,dx_i\right)^{c_i}
\end{eqnarray*}
for any $z\in \R^n$; therefore, we may choose $z_i=S_i(0)$ for $i=1,\ldots,k$ in Lemma~\ref{RBLextproduct}.
\proofbox

\section{$h_{i0}$ is Gaussian in Proposition~\ref{Thetasplits}}
\label{secDependentGaussian}

For positive $C^\alpha$ probability density functions $f$ and $g$
 on $\R^n$ for $\alpha\in (0,1)$, the $C^1$ Brenier map $T:\R^n\to\R^n$ in Theorem~\ref{Breniermap}
 pushing forward the
 the measure on $\R^n$ induced by $g$ to the measure associated to $f$
satisfies that $\nabla T$ is positive definite.
We deduce that
\begin{equation}
\label{Brenierproperty}
\langle T(y)-T(x),y-x\rangle=\int_0^1\langle \nabla T(x+t(y-x))\cdot(y-x),y-x\rangle\,dt\geq 0\mbox{ \ \ for any $x,y\in\R^n$}.
\end{equation}

We say that a continuous function $T:\R^n\to\R^m$ has linear growth if
there exists a positive constant $c>0$ such that
$$
\|T(x)\|\leq c\sqrt{1+\|x\|^2}
$$
for $x\in\R^n$. It is equivalent saying that
\begin{equation}
\label{Tlineargrowthdef}
\limsup_{\|x\|\to\infty}\frac{\|T(x)\|}{\|x\|}<\infty.
\end{equation}
In general, $T$ has polynomial growth, if there exists $k\geq 1$ such that
$$
\limsup_{\|x\|\to\infty}\frac{\|T(x)\|}{\|x\|^k}<\infty.
$$

Proposition~\ref{Tlineargrowth} related to Caffarelli Contraction Principle in  Caffarelli \cite{Caf00}
 was proved by Emanuel Milman,
see for example Colombo, Fathi \cite{CoF21}, De Philippis, Figalli \cite{DPF17},
Fathi, Gozlan, Prod'homme \cite{FGP20}, Y.-H. Kim, E. Milman \cite{KiM12},
Klartag, Putterman \cite{KlP}, Kolesnikov \cite{Kol13}, Livshyts \cite{Liv21} for relevant results.

\begin{prop}[Emanuel Milman]
\label{Tlineargrowth}
If a Gaussian probability density $g$ and a positive $C^\alpha$, $\alpha\in(0,1)$, probability density $f$ on $\R^n$
satisfy $f\leq c\cdot g$ for some positive constant $c>0$,
then the Brenier map $T:\R^n\to\R^n$ pushing forward the measure on $\R^n$ induced by $g$ to the measure associated to $f$
has linear growth.
\end{prop}
\proof We may assume that $g(x)=e^{-\pi\|x\|^2}$.

We observe that $T:\R^n\to\R^n$ is bijective as both $f$ and $g$ are positive. Let $S$ be the inverse of $T$; namely, $S:\R^n\to\R^n$ is the bijective Brenier map
 pushing forward the measure on $\R^n$ induced by $f$ to the measure associated to $g$. In particular,
any Borel $X\subset \R^n$ satisfies
\begin{equation}
\label{BrenierS}
\int_{S(X)}g=\int_X f.
\end{equation}
We note that \eqref{Tlineargrowthdef}, and hence Proposition~\ref{Tlineargrowth} is equivalent saying that
\begin{equation}
\label{TlineargrowS}
\liminf_{x\to\infty}\frac{\|S(x)\|}{\|x\|}>0.
\end{equation}

The main idea of the argument is the following observation.
For any unit vector $u$ and $\theta\in(0,\pi)$, we consider
$$
\Xi(u,\theta)=\left\{y:\langle y,u\rangle \geq \|y\|\cdot \cos\theta\right\}.
$$
Since $S$ is surjective, and $\langle S(z)-S(w), z-w\rangle\geq 0$ for any $z,w\in\R^n$
according to \eqref{Brenierproperty}, we deduce that
 \begin{equation}
\label{SBreniercond}
S(w)+\Xi(u,\theta)\subset S\left(w+\Xi\left(u,\theta+\frac{\pi}2\right)\right)
\end{equation}
for any $u\in S^{n-1}$ and $\theta\in(0,\frac{\pi}2)$.

We suppose that $T$ does not have linear growth, and seek a contradiction. According to
\eqref{TlineargrowS}, there exists a sequence $\{x_k\}$ of points of $\R^n\backslash\{0\}$ tending to infinity such that
$$
\lim_{k\to\infty}\|x_k\|=\infty\mbox{ \ and \ }\lim_{k\to\infty}\frac{\|S(x_k)\|}{\|x_k\|}=0.
$$
In particular, we may assume that
 \begin{equation}
\label{Tlineargrowthno}
\|S(x_k)\|<\frac{\|x_k\|}8.
\end{equation}

For any $k$, we consider the unit vector $e_k=x_k/\|x_k\|$. We observe that
$X_k=x_k+\Xi(e_k,\frac{3\pi}4)$ avoids the interior of the ball $\frac{\|x_k\|}{\sqrt{2}}B^n$; therefore,
if $k$ is large, then
\begin{equation}
\label{fXk}
\int_{X_k} f\leq c\cdot n\kappa_n\int_{\|x_k\|/\sqrt{2}}^\infty r^{n-1}e^{-\pi r^2}\,dr<
\int_{\|x_k\|/\sqrt{2}}^\infty e^{-2 r^2}\sqrt{2} r\,dr
=e^{-\|x_k\|^2}
\end{equation}

On the other hand, $S(x_k)+\Xi(e_k,\frac{\pi}4)$ contains the ball
$$
\widetilde{B}=S(x_k)+\frac{x_k}8+\frac{\|x_k\|}{8\sqrt{2}}\,B^n\subset \frac{\|x_k\|}{2}\,B^n
$$
where we have used \eqref{Tlineargrowthno}. It follows form \eqref{BrenierS} and \eqref{SBreniercond} that
if $k$ is large, then
$$
\int_{X_k} f=\int_{S(X_k)}g \geq \int_{\widetilde{B}}g \geq \kappa_n\left(\frac{\|x_k\|}{8\sqrt{2}}\right)^n
e^{-\pi(\|x_k\|/2)^2}>e^{-\|x_k\|^2}.
$$
This inequality contradicts \eqref{fXk}, and in turn proves \eqref{TlineargrowS}.
\proofbox

Proposition~\ref{dependentGaussian} shows that if the whole space is the dependent subspace and the
Brenier maps corresponding to the
extremizers $f_1,\ldots,f_k$ in Proposition~\ref{Thetasplits} have at most linear growth, then
each $f_i$ is actually Gaussian.
The proof of Proposition~\ref{dependentGaussian} uses classical Fourier analysis, and we refer to Grafakos \cite{Gra14} for the main properties. For our purposes, we need only the action of a tempered distribution on
the space of $C_0^\infty(\R^m)$ of $C^\infty$ functions with compact support, do not need to consider the space of Schwarz functions in general.
We recall that if $u$ is a tempered distribution on Schwarz functions on $\R^n$, then the support
${\rm supp}\,u$ is the intersection of all closed sets $K$ such that if $\varphi\in C_0^\infty(\R^n)$
with ${\rm supp}\,\varphi\subset \R^n\backslash K$, then $\langle u,\varphi\rangle=0$. We write $\hat{u}$ to denote the Fourier transform of a $u$. In particular, if $\theta$ is a function of polynomial growth and $\varphi\in C_0^\infty(\R^n)$, then
\begin{equation}
\label{FourierTransform}
\langle \hat{\theta},\varphi\rangle=\int_{\R^n}\int_{\R^n}\theta(x)\varphi(y)e^{-2\pi i\langle x,y\rangle}\,dxdy.
\end{equation}
We consider the two well-known statements Lemma~\ref{suppsubspace} and Lemma~\ref{supppoint}
 about the support of a Fourier transform to prepare the proof of Proposition~\ref{dependentGaussian}.

\begin{lemma}
\label{suppsubspace}
If $\theta$ is a measurable function of polynomial growth on $\R^n$, and there exist
linear subspace $E$ with $1\leq{\rm dim}\,E\leq n-1$ and function $\omega$ on $E$ such that
$\theta(x)=\omega(P_Ex)$, then
${\rm supp}\, \hat{\theta}\subset E$.
\end{lemma}
\proof We write a $z\in \R^n$ in the form $z=(z_1,z_2)$ with $z_1\in E$ and $z_2\in E^\bot$.
Let $\varphi\in C_0^\infty(\R^n)$ satisfy that ${\rm supp}\,\varphi\subset \R^n\backslash E$, and hence
$\varphi(x_1,o)=0$ for $x_1\in E$, and
the Fourier Integral Theorem in $E^\bot$ implies
$$
\varphi(x_1,z)=\int_{E^\bot}\int_{E^\bot}\varphi(x_1,x_2) e^{2\pi i\langle z-x_2,y_2\rangle}\,dx_2dy_2
$$
for $x_1\in E$ and $z\in E^\bot$.
It follows from \eqref{FourierTransform} that
\begin{eqnarray*}
\langle \hat{\theta},\varphi\rangle&=&
\int_{E^\bot}\int_{E}\int_{E^\bot}\int_{E}\omega(x_1)\varphi(x_1,x_2)e^{-2\pi i\langle x_1,y_1\rangle}
e^{-2\pi i\langle x_2,y_2\rangle}\,dx_1dx_2dy_1dy_2\\
&=&
\int_{E}\int_{E}\omega(x_1)e^{-2\pi i\langle x_1,y_1\rangle}
\left(\int_{E^\bot}\int_{E^\bot}\varphi(x_1,x_2) e^{2\pi i\langle -x_2,y_2\rangle}\,dx_2dy_2\right)dy_1dx_1\\
&=&
\int_{E}\int_{E}\omega(x_1)e^{-2\pi i\langle x_1,y_1\rangle}
\varphi(x_1,0)\,dy_1dx_1=0. \proofbox
\end{eqnarray*}

Next, Lemma~\ref{supppoint} directly follows from Proposition~2.4.1 in Grafakos \cite{Gra14}.

\begin{lemma}
\label{supppoint}
If $\theta$ is a continuous function of polynomial growth on $\R^n$ and
${\rm supp}\, \hat{\theta}\subset\{0\}$, then $\theta$ is a polynomial.
\end{lemma}

\begin{prop}
\label{dependentGaussian}
For  linear subspaces  $E_1,\ldots,E_k$ of $\R^m$ and $c_1,\ldots,c_k>0$ satisfying
\eqref{highdimcond0}, we assume that
\begin{equation}
\label{wholedependent}
\cap_{i=1}^k(E_i\cup E_i^\bot)=\{0\}.
\end{equation}
Let $g_i(x)=e^{-\pi\|x\|^2}$ for $i=1,\ldots,k$ and $x\in E_i$, let
 equality hold in \eqref{RBL}
for positive $C^1$ probability densities $f_i$ on $E_i$, $i=1,\ldots,k$, and
let $T_i:E_i\to E_i$ be the $C^2$ Brenier map satisfying
\begin{equation}
\label{dependentGaussianBrenier}
g_i(x)=\det \nabla T_i(x) \cdot f_i(T_i(x))\mbox{ \ for all $x\in E_i$}.
\end{equation}
If each $T_i$, $i=1,\ldots,k$, has linear growth,
then there exist a positive definite matrix $A:\R^m\to \R^m$ whose eigenspaces are critical subspaces, and
$a_i>0$ and $b_i\in E_i$, $i=1,\ldots,k$, such that
$$
f_i(x)=a_ie^{-\langle Ax,x+b_i\rangle}\mbox{ \ for $x\in E_i$.}
$$
 \end{prop}
\proof We may assume that each linear subspace is non-zero.

We note that the condition \eqref{wholedependent} is equivalent saying that $\R^m$ itself is the dependent subspace with respect to the Brascamp-Lieb data. We may assume that for some $1\leq l\leq k$, we have $1\leq {\rm dim}E_i\leq m-1$
if $i=1,\ldots,l$, and still
\begin{equation}
\label{wholedependent0}
\cap_{i=1}^l(E_i\cup E_i^\bot)=\{0\}.
\end{equation}

We use the diffeomorphism $\Theta:\R^m\to\R^m$ of Proposition~\ref{Thetasplits} defined by
$$
\Theta(y)=\sum_{i=1}^kc_iT_i\left(P_{E_i}y\right),\qquad y\in\R^m.
$$
It follows from \eqref{RBLtranslation} that we may asssume
\begin{equation}
\label{ThetaTioo}
T_i(0)=0\mbox{ \ for $i=1,\ldots,k$, and hence }\Theta(0)=0.
\end{equation}

We claim that there exists a positive definite matrix $B:\R^m\to \R^m$ whose eigenspaces are critical subspaces, and
\begin{equation}
\label{ThetaBclaim}
\nabla\Theta(y)=B\mbox{ \ for $y\in\R^m$}.
\end{equation}
Let $\Theta(y)=(\theta_1(y),\ldots,\theta_m(y))$ for $y\in\R^m$ and $\theta_j\in C^2(\R^m)$, $j=1,\ldots,m$. Since each $T_i$, $i=1,\ldots,k$ has linear growth, it follows that $\Theta$ has linear growth, and in turn each $\theta_j$, $j=1,\ldots,m$, has linear growth.

According to Proposition~\ref{Thetasplits} (iii),
 there exist $C^2$ functions $\Omega_i:E_i\to E_i$
and $\Gamma_i:E_i^\bot\to E_i^\bot$ such that
$$
\Theta(y)=\Omega_i(P_{E_i}y)+\Gamma_i(P_{E_i^\bot}y)
$$
for $i=1,\ldots,k$ and $y\in \R^n$. We write $\Omega_i(x)=(\omega_{i1}(x),\ldots,\omega_{im}(x))$
and $\Gamma_i(x)=(\gamma_{i1}(x),\ldots,\gamma_{im}(x))$; therefore,
\begin{equation}
\label{thetajslip}
\theta_j(y)=\omega_{ij}(P_{E_i}y)+\gamma_{ij}(P_{E_i^\bot}y)
\end{equation}
for $j=1,\ldots,m$ and $i=1,\ldots,k$.

Fix a $j\in\{1,\ldots,m\}$. It follows from Lemma~\ref{suppsubspace} and \eqref{thetajslip} that
$$
{\rm supp}\,\hat{\theta}_j\subset E_i\cup E_i^\bot
$$
for $i=1,\ldots,l$. Thus \eqref{wholedependent0} yields that
$$
{\rm supp}\,\hat{\theta}_j\subset \{0\},
$$
and in turn we deduce from Lemma~\ref{supppoint} that $\theta_j$ is a polynomial. Given that $\theta_j$ has linear growth, it follows that there exist $w_j\in \R^m$ and $\alpha_j\in\R$ such that $\theta_j(y)=\langle w_j,y\rangle+\alpha_j$.
We deduce from $\theta_j(o)=0$ ({\rm cf.} \eqref{ThetaTioo}) that $\alpha_j=0$.

The argument so far yields that there exists an $m\times m$ matrix $B$ such that $\Theta(y)=By$ for $y\in\R^m$.
As $\nabla  \Theta(y)=B$ is positive definite and its eigenspaces are critical subspaces, we conclude
the claim \eqref{ThetaBclaim}.

Since
$\nabla T_i(P_{E_i}y)=\nabla\Theta(y)|_{E_i}$ for $i=1,\ldots,k$ and $y\in\R^m$ by Proposition~\ref{Thetasplits} (iv), we deduce that $T_i^{-1}=B^{-1}|_{E_i}$ for $i=1,\ldots,k$. It follows from \eqref{dependentGaussianBrenier} that
$$
f_i(x)= e^{-\pi\|B^{-1}x\|^2}\cdot \det \left(B^{-1}|_{E_i}\right)\mbox{ \ for $x\in E_i$}
$$
for $i=1,\ldots,k$. Therefore, we can choose $A=\pi B^{-2}$.
\proofbox

\section{Proof of Theorem~\ref{RBLtheoequa} }
\label{secFinalProof}

We may assume that each linear subspace $E_i$ is non-zero in Theorem~\ref{RBLtheoequa}. Let $f_i$ be a probability density on $E_i$ in a way such that
equality holds for $f_1,\ldots,f_k$ in \eqref{RBL}. For $i=1,\ldots,k$ and $x\in E_i$, let $g_i(x)=e^{-\pi\|x\|^2}$, and hence
 $g_i$
is a probability distribution on $E_i$, and
$g_1,\ldots,g_k$ are extremizers in  Barthe's  inequality \eqref{RBL}.

It follows from Lemma~\ref{convolution} that the convolutions $f_1*g_1,\ldots,f_k*g_k$ are also extremizers for \eqref{RBL}. We observe that for $i=1,\ldots,k$, $f_i*g_i$ is a bounded positive $C^\infty$ probability density on $E_i$. Next we deduce from
Lemma~\ref{RBLextproduct} that there exist $z_i\in E_i$ and $\gamma_i>0$ for $i=1,\ldots,k$ such that defining
$$
\tilde{f}_i(x)=\gamma_i\cdot g_i(x)\cdot (f_i*g_i)(x-z_i) \mbox{ \ for $x\in E_i$},
$$
$\tilde{f}_1,\ldots,\tilde{f}_k$ are probability densities that are extremizers for \eqref{RBL}. We note that if $i=1,\ldots,k$, then
$\tilde{f}_i$ is positive and $C^\infty$, and there exists $c>1$ satisfying
\begin{equation}
\label{tildefcg}
\tilde{f}_i\leq c\cdot g_i.
\end{equation}

Let $\widetilde{T}_i:E_i\to E_i$ be the $C^\infty$ Brenier map satisfying
\begin{equation}
\label{MongeAmperegtildef}
g_i(x)=\det \nabla \widetilde{T}_i(x) \cdot \tilde{f}_i(\widetilde{T}_i(x))\mbox{ \ for all $x\in E_i$},
\end{equation}
We deduce from \eqref{tildefcg} and Proposition~\ref{Tlineargrowth} that $\widetilde{T}_i$ has linear growth.

For $i=1,\ldots,k$ and $x\in F_0\cap E_i$, let $g_{i0}(x)=e^{-\pi\|x\|^2}$.
It follows from Proposition~\ref{Thetasplits} (i) that
 for $i\in\{1,\ldots,k\}$, there exists positive $C^1$
integrable $h_{i0}:\,F_0\cap E_i\to[0,\infty)$ (where $h_{i0}(o)=1$ if $F_0\cap E_i=\{0\}$),
and for any $i\in\{1,\ldots,k\}$ and $j\in\{1,\ldots,l\}$ with $F_j\subset E_i$,
 there exists positive $C^1$
integrable $\tilde{h}_{ij}:\,F_j\to[0,\infty)$  such that
$$
\tilde{f}_i(x)=\tilde{h}_{i0}(P_{F_0}x)\cdot
\prod_{F_j\subset E_i\atop j\geq 1}\tilde{h}_{ij}(P_{F_j}x) \mbox{ \ \ \  for $x\in E_i$}.
$$
We deduce from Proposition~\ref{Thetasplits} (ii) that $\widetilde{T}_{i0}=\widetilde{T}_i|_{F_0\cap E_i}$ is the Brenier map pushing forward
the measure on $F_0\cap E_i$ determined $g_{i0}$ onto the measure determined by $\tilde{h}_{i0}$. Since
$\widetilde{T}_i$ has linear growth, $\widetilde{T}_{i0}$ has linear growth, as well, for $i=1,\ldots,k$.

We deduce from Proposition~\ref{Thetasplitshij} (i)
that $\sum_{i=1}^kc_iP_{E_i\cap F_0}={\rm Id}_{F_0}$,
the Geometric Brascamp Lieb data $E_1\cap F_0,\ldots,E_k\cap F_0$ in $F_0$ has no independent subspaces, and
$\tilde{h}_{10},\ldots,\tilde{h}_{k0}$ are extremizers in  Barthe's  inequality for this data in $F_0$.

As $\widetilde{T}_{i0}$ has linear growth for $i=1,\ldots,k$,
Proposition~\ref{dependentGaussian} yields the existence of a positive definite matrix
$\widetilde{A}:F_0\to F_0$ whose eigenspaces are critical subspaces, and
$\tilde{a}_i>0$ and $\tilde{b}_i\in F_0\cap E_i$ for $i=1,\ldots,k$, such that
$$
\tilde{f}_i(x)=\tilde{a}_ie^{-\langle \widetilde{A}x,x+\tilde{b}_i\rangle}\cdot
\prod_{F_j\subset E_i\atop j\geq 1}\tilde{h}_{ij}(P_{F_j}x) \mbox{ \ \ \  for $x\in E_i$}.
$$
Dividing by $g_i$ and shifting, we deduce that there exist a symmetric
matrix
$\bar{A}:F_0\to F_0$ whose eigenspaces are critical subspaces, and
$\bar{a}_i>0$ and $\bar{b}_i\in F_0\cap E_i$ for $i=1,\ldots,k$,
and for any $i\in\{1,\ldots,k\}$ and $j\in\{1,\ldots,l\}$ with $F_j\subset E_i$,
 there exists positive $C^1$
 $\bar{h}_{ij}:\,F_j\to[0,\infty)$  such that
$$
f_i*g_i(x)=\bar{a}_ie^{-\langle \bar{A}x,x+\bar{b}_i\rangle}\cdot
\prod_{F_j\subset E_i\atop j\geq 1}\bar{h}_{ij}(P_{F_j}x) \mbox{ \ \ \  for $x\in E_i$}.
$$
Since $f_i*g_i$ is a probability density on $E_i$, it follows that
$\bar{A}$ is positive definite and $\bar{h}_{ij}\in L_1(E_i\cap F_j)$  for
$i\in\{1,\ldots,k\}$ and $j\in\{1,\ldots,l\}$ with $F_j\subset E_i$.

 For any $i=1,\ldots,k$, we write $\hat{\varrho}$ for the Fourier transform of a function $\varrho\in L_1(E_i)$, thus we can take the inverse Fourier transform in the sense that
$\varrho$ is a.e. the $L_1$ limit of
$$
x\mapsto \int_{\R^n}\hat{\varrho}(\xi)e^{-a|\xi|^2}e^{2\pi i\langle \xi, x\rangle}\,d\xi
$$
as $a>0$ tends to zero.
For $i=1,\ldots,k$, using that $\widehat{f_i*g_i}=\hat{f}_i\cdot \hat{g}_i$, we deduce that the restriction of $\hat{f}_i$ to $F_0\cap E_i$ is the quotient of two Gaussian densities. Since $\hat{f}_i$ is bounded and zero at infinity, we deduce that
the restriction of $\hat{f}_i$ to $F_0\cap E_i$ is a Gaussian density
for $i=1,\ldots,k$, as well, with the symmetric matrix involved being positive definite.
We conclude using the inverse Fourier transform above and the fact that the linear subspaces $F_j$, $j=0,\dots,l$, are pairwise orthogonal that
there exist a symmetric
matrix
$A:F_0\to F_0$ whose eigenspaces are critical subspaces, and
$a_i>0$ and $b_i\in F_0\cap E_i$ for $i=1,\ldots,k$,
and for any $i\in\{1,\ldots,k\}$ and $j\in\{1,\ldots,l\}$ with $F_j\subset E_i$,
 there exists
 $h_{ij}:\,F_j\to[0,\infty)$  such that
$$
f_i(x)=a_ie^{-\langle Ax,x+b_i\rangle}\cdot
\prod_{F_j\subset E_i\atop j\geq 1}h_{ij}(P_{F_j}x) \mbox{ \ \ \  for a.e. $x\in E_i$}.
$$
Since $f_i$ is a probability density on $E_i$, it follows that $A$ is positive definite and
each $h_{ij}$ is non-negative and integrable.
Finally,
Proposition~\ref{Thetasplitshij} (ii) yields that
there exist integrable $\psi_{j}:\,F_j\to[0,\infty)$
for $j=1,\ldots,l$ where $\psi_j$ is log-concave whenever $F_j\subset E_\alpha\cap E_\beta$ for $\alpha\neq \beta$, and there exist $a_{ij}>0$ and $b_{ij}\in F_j$
for any $i\in\{1,\ldots,k\}$ and $j\in\{1,\ldots,l\}$ with $F_j\subset E_i$ such that
$h_{ij}(x)=a_{ij}\cdot \psi_j(x-b_{ij})$ for $i\in\{1,\ldots,k\}$ and $j\in\{1,\ldots,l\}$ with $F_j\subset E_i$.

Finally, we assume that $f_1,\ldots,f_k$ are of the form as described in \eqref{RBLtheoequaform} and equality holds for  all
$x\in E_i$ in \eqref{RBLtheoequaform}. According to \eqref{RBLtranslation}, we may assume that there
exist a positive definite matrix $\Phi:F_0\to F_0$ whose proper eigenspaces are critical subspaces
and a $\tilde{\theta}_i>0$ for $i=1,\ldots,k$ such that
\begin{equation}
\label{RBLtheoequaform0}
f_i(x)=\tilde{\theta}_i e^{-\|\Phi P_{F_0}x\|^2}\prod_{F_j\subset E_i}h_{j}(P_{F_j}(x)) \mbox{ \ \ \  for $x\in E_i$}.
\end{equation}
We recall that according to \eqref{sumci1}, if $j\in\{1,\ldots,l\}$, then
\begin{equation}
\label{sumci10}
\sum_{E_i\supset F_j}c_i=1.
\end{equation}
We set $\theta=\prod_{i=1}^k\tilde{\theta}_i^{c_i}$ and $h_0(x)=e^{-\|\Phi x\|^2}$  for   $x\in F_0$.
On the left hand side of  Barthe's  inequality \eqref{RBL}, we
  use first \eqref{sumci10} and the log-concavity of
$h_j$ whenever $j\geq 1$ and $F_j\subset E_\alpha\cap E_\beta$ for $\alpha\neq \beta$,
 secondly Proposition~\ref{GaussianExtremizers}, thirdly
\eqref{sumci10}, fourth  the Fubini Theorem, and finally \eqref{sumci10} again to prove that
\begin{eqnarray*}
\int_{*,\R^n}\sup_{x=\sum_{i=1}^kc_ix_i\atop x_i\in E_i}\;\prod_{i=1}^kf_i(x_i)^{c_i}\,dx&=&
\theta\int_{*,\R^n}\sup_{x=\sum_{i=1}^k\sum_{j=0}^lc_ix_{ij}\atop x_{ij}\in E_i\cap F_j}\;
\prod_{j=0}^l\prod_{i=1}^kh_j(x_{ij})^{c_i}\,dx\\
&=& \theta\int_{*,\R^n}\prod_{j=0}^l\sup_{P_{F_j}x=\sum_{i=1}^kc_ix_{ij}\atop x_{ij}\in E_i\cap F_j}\;
\prod_{i=1}^kh_j(x_{ij})^{c_i}\,dx\\
&=& \theta\int_{*,\R^n}
\left(\sup_{P_{F_0}x=\sum_{i=1}^kc_ix_{i0}\atop x_{i0}\in E_i\cap F_0}
\prod_{i=1}^k e^{-c_i\|\Phi x_{i0}\|^2}\right)\times
\prod_{j=1}^lh_j(P_{F_j}x)\,dx\\
&=&\theta\left(\prod_{i=1}^k\left(\int_{F_0\cap E_i}e^{-\|\Phi y\|^2}\,dy\right)^{c_i}\right)
\times \prod_{j=1}^l\int_{F_j}h_j\\
&=&\prod_{i=1}^k\left(\int_{E_i}f_i\right)^{c_i},
\end{eqnarray*}
completing the proof of Theorem~\ref{RBLtheoequa}.
\proofbox

\section{Equality in the Bollobas-Thomason inequality and in its dual}
\label{secBT}

We fix an orthonormal basis $e_1,\ldots,e_n$ of $\R^n$ for the whole section. We set
$\sigma_i^0=\sigma_i$ and $\sigma_i^1=[n]\setminus\sigma_i$. When we write
$\tilde{\sigma_1},\ldots,\tilde{\sigma_l}$ for the induced cover from $\sigma_1,\ldots,\sigma_k$, we assume that the sets
$\tilde{\sigma_1},\ldots,\tilde{\sigma_l}$ are pairwise distinct.

\begin{lemma}\label{combdesc}
 For $s\geq 1$, let
$\sigma_1,\ldots,\sigma_k\subset[n]$ form an $s$-uniform cover of $[n]$, and let
$\tilde{\sigma}_1,\ldots,\tilde{\sigma}_\ell$ be
the $1$-uniform cover  of $[n]$ induced by
 $\sigma_1,\ldots,\sigma_k$. Then
\begin{description}
\item{(i)} the subspaces $E_{\sigma_i}:={\rm lin}\{e_j:i\in\sigma_i\}$ satisfy
\begin{align}\label{GBLdata}
\displaystyle \sum_{i=1}^k \frac1s\,P_{E_{\sigma_i}}=I_n
\end{align}
i.e. form a Geometric Brascamp Lieb data;
\item{(ii)} For $r\in\tilde{\sigma}_j$, $j=1,\ldots,\ell$, we have
\begin{align}
\label{jhjbvsvjs}
\tilde{\sigma}_j:=\bigcap_{r\in\sigma_i}\sigma_i^0\cap\bigcap_{r\notin\sigma_i}\sigma_i^1;
\end{align}
\item{(iii)}
the subspaces $F_{\tilde{\sigma}_j}:={\rm lin}\{e_r:r\in\tilde{\sigma}_j\}$ are the independent subspaces of the Geometric Brascamp Lieb data \eqref{GBLdata} and
$F_{\rm dep}=\{0\}$.
\end{description}
\end{lemma}
\proof  Since $\sigma_1,\ldots,\sigma_k$ form a $s$-uniform cover, every $e_i\in\R^n$ is contained in exactly $s$ of $E_{\sigma_1},\ldots,E_{\sigma_k}$, yielding (i).

For (ii), the definition of $\tilde{\sigma}_j$ directly implies \eqref{jhjbvsvjs}.

For (iii), the linear subspaces $F_{\tilde{\sigma}_1},\ldots, F_{\tilde{\sigma}_\ell}$  are pairwise orthogonal because
$\sigma_i^0\cap\sigma_i^1=\emptyset$ for $i=1,\ldots,k$. On the other hand, for any $r\in[n]$,
$r\in\cap_{i=1}^n \sigma^{\varepsilon(i)}_i$ where $\varepsilon(i)=0$ if $r\in\sigma_i$, and
$\varepsilon(i)=1$ if $r\not\in\sigma_i$; therefore, $F_{\tilde{\sigma}_1},\ldots, F_{\tilde{\sigma}_\ell}$
span $\R^n$. In particular, $F_{\rm dep}=\{0\}$.
\hfill \proofbox

Let us introduce the notation that we use when handling both the Bollobas-Thomason inequality and its dual.
Let $\sigma_1,\ldots,\sigma_k$ be the $s$ cover of $[n]$ occuring in
Theorem~\ref{Bollobas-Thomason-eq} and Theorem~\ref{Liakopoulos-eq}, and hence
$E_i=E_{\sigma_i}$, $i=1,\ldots,k$, satisfies
\begin{equation}
\label{sumPEsigmai}
\frac1s\sum_{i=1}^k P_{E_{\sigma_i}}=I_n.
\end{equation}
Let $\tilde{\sigma}_1,\ldots,\tilde{\sigma}_l$ be the $1$-uniform cover of $[n]$ induced by
$\sigma_1,\ldots,\sigma_k$. It follows that
\begin{eqnarray}
\label{independent}
F_j&=&E_{\tilde{\sigma}_j} \mbox{ \ \ for $j=1,\ldots,l$ are the independent subspaces},\\
\label{dependent}
F_{\rm dep}&=&\{0\}.
\end{eqnarray}
For any $i\in\{1,\ldots,k\}$, we set
$$
I_i=\{j\in\{1,\ldots,l\}:\,F_j\subset E_i\},
$$
and for any $j\in\{1,\ldots,l\}$, we set
$$
J_j=\{i\in\{1,\ldots,k\}:\,F_j\subset E_i\}.
$$

For the reader's convenience, we restate Theorem~\ref{Bollobas-Thomason} and
Theorem~\ref{Bollobas-Thomason-eq} as Theorem~\ref{Bollobas-Thomason1},
and Theorem~\ref{Liakopoulos} and Theorem~\ref{Liakopoulos-eq} as Theorem~\ref{Liakopoulos1}.

\begin{theo}
	\label{Bollobas-Thomason1}
	If $K\subset \R^n$ is compact and affinely spans $\R^n$, and
	$\sigma_1,\ldots,\sigma_k\subset[n]$ form an $s$-uniform cover of $[n]$ for $s\geq 1$, then
	\begin{equation}
	\label{Bollobas-Thomasson-ineq1}
	|K|^s\leq \prod_{i=1}^k|P_{E_{\sigma_i}}K|.
	\end{equation}
	Equality holds if and only if
	$K=\oplus_{i=1}^l P_{F_{\tilde{\sigma}_i}}K$
	where $\tilde{\sigma}_1,\ldots,\tilde{\sigma}_l$ is
	the $1$-uniform cover  of $[n]$ induced by
	$\sigma_1,\ldots,\sigma_k$ and $F_{\tilde{\sigma_i}}$ is the linear hull of the $e_i$'s with indeces from $\tilde{\sigma_i}$.
\end{theo}
\proof
We set $E_i:=E_{\sigma_i}$ which subspaces compose a geometric data according to Lemma \ref{combdesc}. We start with a proof of Bollobas-Thomason inequality. It follows directly
from the Brascamp-Lieb inequality as
\begin{align}
\nonumber
|K|=\int_{\R^n}1_K(x)\,dx &\leq \int_{\R^n}\prod_{i=1}^{k}1_{P_{E_i}(K)}(P_{E_i}(x))^{\frac{1}{s}}\,dx\\
\label{BTproof}
&\leq\prod_{i=1}^{k}\Big(\int_{E_i}1_{P_{E_i}(K)}\Big)^{\frac{1}{s}}=\prod_{i=1}^{k}|P_{E_i}(K)|^{\frac{1}{s}}
\end{align}
where the first inequality is from the monotonicity of the integral while the second is Brascamp-Lieb inequality Theorem~\ref{BLtheo}. Now, if equality holds in \eqref{BTproof}, then on the one hand,
\begin{align*}
1_K(x)=\prod_{i=1}^{k}1_{P_{E_i}(K)}(P_{E_i}(x))
\end{align*}
and on the other hand, if $F_1,\ldots, F_l$ are the independent subspaces of the data, then they  span $\R^n$ according to Lemma~\ref{combdesc}; namely, $F_{{\rm dep}}=\{0\}$. It follows from Theorem~\ref{BLtheoequa} that there exist integrable functions $h_j:F_j\to \R$, such that, for Lebesgue a.e. $x_i\in E_i$
\begin{align*}
1_{P_{E_i}K}(x_i)=\theta_i\prod_{j\in I_i}h_j(P_{F_j}(x_i))
\end{align*}
Therefore from the previous two, we have for $x\in\R^n$
\begin{align*}
1_{K}(x)=\prod_{i=1}^{k}\theta_i\prod_{j\in I_i}h_j(P_{F_j}(P_{E_i}(x)))
\end{align*}
Now, since for $j\in I_i$ we have $F_j\subset E_i$ we can delete the $P_{E_i}$ on the above product. Thus, for $\theta=\prod_{i=1}^k\theta_i$, we have for Lebesgue a.e. $x\in\R^n$
\begin{equation}
\label{1Kprodhjs}
1_K(x)=\theta\prod_{i=1}^k\prod_{j\in I_i}h_j(P_{F_j}(x))
=\theta\prod_{j=1}^l h_j(P_{F_j}(x))^{|J_j|}.
\end{equation}
Now, for $x\in K$ the last product on above is constant, so
\begin{align}\label{theta}
\theta=\frac{1}{\prod_{i=1}^lh_j(P_{F_j}(x_0))^{|J_j|}}
\end{align}
for some $x_o\in K$. For $j=1,\ldots,l$ we set $\varphi_j:F_j\to\R^n$, by
$$
\varphi_j(x)=\frac{h_j(x+P_{F_j}(x_0))^{|J_j|}}{h_j(P_{F_j}(x_0))^{|J_j|}}.
$$
We see that $\varphi_j(o)=1$ and also
\eqref{1Kprodhjs} and \eqref{theta} yields
\begin{equation}
\label{1Kprodphijs}
1_{K-x_0}(x)=\prod_{j=1}^l \varphi_j(P_{F_j}(x))
\end{equation}
For $m\in\{1,\ldots,l\}$,
taking $x\in F_m$ in \eqref{1Kprodphijs} (and hence $\varphi_j(P_{F_j}(x))=1$ for $j\neq m$)
shows that
$$
1_{K-x_0}(y)=\varphi_m(y),
$$
for Lebesgue a.e. $y\in F_m$. Therefore \eqref{1Kprodphijs} and the ortgonality of the $F_j$'s,
$$
K-x_0=\bigcap_{j=1}^lP_{F_j}^{-1}(P_{F_j}(K-x_o))=\bigoplus_{j=1}^lP_{F_j}(K-x_o),
$$
completing the proof of Theorem~\ref{Bollobas-Thomason1}.
\hfill \proofbox

To prove Theorem~\ref{Liakopoulos1}, we use two small observations. First
if $M$ is any convex body with $o\in{\rm int}\,M$, then
\begin{equation}
\label{VKexponential}
\int_{\R^n}e^{-\|x\|_M}\,dx=\int_0	^\infty e^{-r}nr^{n-1}|M|\,dr=n! |M|.
\end{equation}
Secondly, if $F_j$ are pairwise orthogonal subspaces and $M={\rm conv}\,\{M_1,\ldots,M_l\}$ where $M_j\subset F_j$ is a ${\rm dim}F_j$-dimensional compact convex set with $o\in{\rm relint}\,M_j$, then for any $x\in\R^n$
\begin{equation}
\label{sumMjnorm}
\|x\|_M=\sum_{i=1}^l\|P_{F_j}x\|_{M_j}.
\end{equation}
In addition, we often use the fact, for a subspace $F$ of $\R^n$ and $x\in F$, then $\|x\|_K=\|x\|_{K\cap F}$.

\begin{theo}
	\label{Liakopoulos1}
	If $K\subset \R^n$ is compact convex with $o\in{\rm int}K$, and
	$\sigma_1,\ldots,\sigma_k\subset[n]$ form an $s$-uniform cover of $[n]$ for $s\geq 1$, then
	\begin{equation}
	\label{Liakopoulos-ineq1}
	|K|^s\geq \frac{\prod_{i=1}^k|\sigma_i|!}{(n!)^s}\cdot \prod_{i=1}^k|K\cap E_{\sigma_i}|.
	\end{equation}
	Equality holds if and only if
	$K={\rm conv}\{E_{\tilde{\sigma}_i}\cap K\}_{i=1}^l$
	where $\tilde{\sigma}_1,\ldots,\tilde{\sigma}_l$ is
	the $1$-uniform cover  of $[n]$ induced by
	$\sigma_1,\ldots,\sigma_k$.
\end{theo}
\proof
We define
\begin{align}\label{defoff}
f(x)=e^{-\|x\|_K},
\end{align}
which is a log-concave function with $f(o)=1$, and satisfying ({\it cf} \eqref{VKexponential})
\begin{equation}
\label{ftonint}
\int_{\R^n}f(y)^n\,dy=\int_{R^n}e^{-n\|y\|_K}\,dy=\int_{R^n}e^{-\|y\|_{\frac{1}{n}K}}
=n!\left|\frac{1}{n}K\right|=\frac{n!}{n^n}\cdot |K|.
\end{equation}
We claim that
\begin{align}\label{weaker-RBL}
n^n\int_{\R^n}f(y)^n\,dy\geq\prod_{i=1}^k\Big(\int_{E_i}f(x_i)\,dx_i\Big)^{1/s}.
\end{align}
Equating the traces of the two sides of \eqref{GBLdata}, we deduce
that, $d_i:=|\sigma_i|={\rm dim}E_i$
\begin{equation}
\label{sumdin}
\sum_{i=1}^k\frac{d_i}{sn}=1.
\end{equation}
For $z=\sum_{i=1}^k\frac{1}{s}x_i$ with $x_i\in E_i$, the log-concavity of $f$ and its definition \eqref{defoff}, imply
\begin{equation}
\label{znsupf}
f(z/n)\geq \prod_{i=1}^kf(x_i/d_i)^{\frac{d_i}{ns}}= \prod_{i=1}^kf(x_i)^{\frac{1}{ns}}.
\end{equation}
Now, the monotonicity of the integral and Barthe's inequality yield
\begin{equation}
\label{sumPEsigmai0}
\int_{\R^n}f(z/n)^{n}\,dz\geq
\int_{\R^n}^*\sup_{z=\sum_{i=1}^k\frac{1}{s}x_i,\, x_i\in E_i}\prod_{i=1}^kf(x_i)^{1/s}\,dz
\geq\prod_{i=1}^k\Big(\int_{E_i}f(x_i)\,dx_i\Big)^{1/s}.
\end{equation}
Making the change of variable $y=z/n$ we conclude to \eqref{weaker-RBL}. Computing the right hand side of \eqref{weaker-RBL},  we have
\begin{align}\label{compright}
\int_{E_i}f(x_i)\,dx_i=\int_{E_i}e^{-\|x_i\|_K}\,dx_i=\int_{E_i}e^{-\|x_i\|_{K\cap E_i}}\,dx_i=d_i!|K\cap E_i|.
\end{align}
Therefore, \eqref{ftonint}, \eqref{weaker-RBL} and \eqref{compright} yield \eqref{Liakopoulos-ineq1}.

Let us assume that equality holds in \eqref{Liakopoulos-ineq1}, and hence we have two equalities in \eqref{sumPEsigmai0}. We set
\begin{align*}
M={\rm conv}\{K\cap F_j\}_{1\leq j\leq l}.
\end{align*}
Clearly, $K\supseteq M$. For the other inclusion, we start with $z\in {\rm int}K$, namely $\|z\|_K<1$. Equality in the first inequality in \eqref{sumPEsigmai0} means,
$$
\left(e^{-\|z/n\|_K}\right)^n=\sup_{z=\sum_{i=1}^k\frac{1}{s}x_i,\, x_i\in E_i}\prod_{i=1}^ke^{-\|x_i\|_K1/s},
$$
or in other words,
\begin{equation}
\label{znormKsEi}
\|z\|_K=\frac1s\cdot \inf_{z=\sum_{i=1}^k\frac{1}{s}x_i,\, x_i\in E_i}\sum_{i=1}^k \|x_i\|_K
=\inf_{z=\sum_{i=1}^k y_i,\, y_i\in E_i}\sum_{i=1}^k \|y_i\|_K.
\end{equation}
We deduce that there exist $y_i\in E_i$, $i=1,\ldots,k$ such that
\begin{equation}
\label{zxicond}
z=\sum_{i=1}^k y_i\mbox{ \ and \ }\sum_{i=1}^k \|y_i\|_K<1,
\end{equation}
Therefore, from \eqref{zxicond}, then \eqref{sumMjnorm} and after the triangle inequality for $\|\cdot\|_{K\cap F_j}$, we have
\begin{align}\label{den kserw}
\|z\|_M=\left\|\sum_{i=1}^k\sum_{j\in I_i}P_{F_j}y_i\right\|_{M}=\sum_{i=1}^k\left\|\sum_{i\in I_i}P_{F_j}y_i\right\|_{K\cap F_j}\leq
\sum_{i=1}^k\sum_{i\in I_i}\left\|P_{F_j}y_i\right\|_{K\cap F_j}.
\end{align}
It suffices to show that
\begin{align}\label{cnsldnslvrn}
K\cap E_i={\rm conv}\{K\cap F_j\}_{j\in I_i}
\end{align}
because then, from \eqref{den kserw}, applying \eqref{sumMjnorm} and \eqref{zxicond}, we have
\begin{align*}
\|z\|_M\leq
\sum_{j=1}^l\sum_{i\in J_j}\left\|P_{F_j}y_i\right\|_{K\cap F_j}=\sum_{i=1}^k \|y_i\|_{K\cap E_i}<1,
\end{align*}
which means $z\in M$. Now, to show \eqref{cnsldnslvrn}, we start with the equality case of  Barthe's  inequality which has been applied in \eqref{sumPEsigmai0}. From Theorem~\ref{RBLtheoequa}, there exist
$\theta_i>0$ and $w_i\in E_i$ and
log-concave $h_{j}:\,F_j\to[0,\infty)$, namely $h_j=e^{-\varphi_j}$ for a convex functon $\varphi_j$, such that
\begin{align}\label{csdkcvsrv}
e^{-\|x_i\|_{K\cap E_i}}=\theta_i\prod_{j\in I_i}h_{j}(P_{F_j}(x_i-w_{i})).
\end{align}
for Lebesgue a.e. $x_i\in E_i$. For $i\in[k]$ and $j\in I_i$ we set, $\psi_{ij}:F_j\to \R$ by
\begin{align*}
\psi_{ij}(x)=\varphi_j\left(x-P_{F_j}w_i\right)-\varphi_j\left(-P_{F_j}w_i\right)+\frac{\ln \theta_i}{|I_i|}.
\end{align*}
We see
\begin{align}\label{csdvsjwcsa}
\psi_{ij}(o)=0 \ \text{and} \ \psi_{ij} \ \text{is convex on} \ F_j.
\end{align}
and also \eqref{csdkcvsrv} yields, for $x\in E_i$
\begin{equation}
\label{exKEi}
e^{-\|x\|_{K\cap E_i}}=
\exp\left(-\sum_{j\in I_i}\psi_{ij}(P_{F_j}x)\right).
\end{equation}
For $x\in F_j$, we apply $\lambda x$ to \eqref{exKEi}  with $\lambda>0$, and we have
from $\psi_{im}(o)=0$ for $m\in I_i\backslash \{j\}$ that
\begin{align}\label{key}
\psi_{ij}(\lambda x)=\lambda \psi_{ij}(x)\mbox{ and} \  \psi_{ij}(x)>0.
\end{align}
We deduce from \eqref{csdvsjwcsa} and \eqref{key} that $\psi_{ij}$ is a norm. Therefore,  $\psi_{ij}(x)=\|x\|_{C_{ij}}$ for some $({\rm dim}\,F_j)$-dimensional compact convex set
$C_{ij}\subset F_j$ with $o\in{\rm relint}\,C_{ij}$. Now \eqref{exKEi} becomes,
\begin{align*}
\|x\|_{K\cap E_i}=\sum_{j\in I_i}\|P_{F_j}x\|_{C_{ij}}
\end{align*}
and hence by
\eqref{sumMjnorm} we conclude to
$$
K\cap E_i={\rm conv}\,\{C_{ij}\}_{j\in I_i}.
$$
In particular, if $i\in[k]$ and $j\in I_i$, then
$C_{ij}=(K\cap E_i)\cap F_j=K\cap F_j$,
completing the proof of
 \eqref{cnsldnslvrn}, and in turn yielding Theorem~\ref{Liakopoulos-eq}.
\hfill \proofbox

\noindent{\bf Acknowledgements } We thank Alessio Figalli, Greg Kuperberg and Christos Saroglou for helpful discussions.
We are especially grateful to Emanuel Milman for providing the proof of Proposition~\ref{Tlineargrowth}, and for Franck Barthe for providing the proof of Proposition~\ref{GBLconstanthighdim} and insight on the history of the subject,
and for further ideas and extremely helpful discussions.
We thank the referee for correcting a mistake in Theorem~\ref{RBLtheoequa} and signicantly improving the presentation of the whole paper.

The first named author is also grateful for the hospitality and excellent working environment provided by University of California, Davis and by ETH Z\"urich during various parts of this project.

\end{document}